\documentclass[review]{elsarticle}  
\usepackage{tikz}
\usepackage{amsmath,bm,algorithm,algorithmic}
\usepackage[mathscr]{euscript}
 \let\mathscr\relax
\usepackage{amssymb}
\usepackage{mathrsfs}
\usepackage{graphicx} 
\usepackage{natbib}
\usepackage{lineno,hyperref}
\usepackage{xcolor}
\usepackage{caption}
\usepackage{subcaption}
\usepackage{booktabs}

\begin{document}

\begin{frontmatter}

\title{An efficient adaptive sparse grid collocation method through derivative estimation}

\author{Anindya Bhaduri}
\author{Lori Graham-Brady\footnote{Corresponding author. Email address: lori@jhu.edu}}
\address{Department of Civil Engineering,Johns Hopkins University,3400 N. Charles Street, Baltimore, MD 21218, United Sates}




\begin{abstract}
For uncertainty propagation of highly complex and/or nonlinear problems, one must resort to sample-based non-intrusive approaches \cite{le2010spectral}. In such cases, minimizing the number of function evaluations required to evaluate the response surface is of paramount importance. Sparse grid approaches have proven effective in reducing the number of sample evaluations. For example, the discrete projection collocation method has the notable feature of exhibiting fast convergence rates when approximating smooth functions; however, it lacks the ability to accurately and efficiently track response functions that exhibit fluctuations, abrupt changes or discontinuities in very localized regions of the input domain. On the other hand, the piecewise linear collocation interpolation approach can track these localized variations in the response surface efficiently, but it converges slowly in the smooth regions. The proposed methodology, building on an existing work on adaptive hierarchical sparse grid collocation algorithm \cite{ma2009adaptive}, is able to track localized behavior while also avoiding unnecessary function evaluations in smoother regions of the stochastic space by using a finite difference based one-dimensional derivative evaluation technique in all the dimensions. This derivative evaluation technique leads to faster convergence in the smoother regions than what is achieved in the existing collocation interpolation approaches. Illustrative examples show that this method is well suited to high-dimensional stochastic problems, and that stochastic elliptic problems with stochastic dimension as high as 100 can be dealt with effectively.
\end{abstract}


\end{frontmatter}

\section{Introduction} Any model consists of input parameters which are inherently random. The uncertainty in the inputs naturally leads to an uncertainty in the output. Thus a single solution for the system using a fixed set of input parameters is not sufficient to describe the system completely.  Thus, given the input uncertainties, it is of real interest to understand how these uncertainties propagate through the deterministic system model and result in uncertainties in the output solution. The quantification of the output uncertainties is a much more comprehensive descriptor of the system under study. \\
\indent The traditional approach is to use random sampling techniques such as Monte Carlo (MC) method. It involves generating sets of realizations of all the input parameters following their individual probability distributions and then solving the deterministic code for each set of realizations. The advantage of this method is that it is easy to implement, it has a non-intrusive nature and the convergence rate is independent of the number of stochastic dimensions. On the other hand, it suffers from the drawback that it cannot easily approximate the solution space and usually only gives the output statistics, such as the mean and the variance. The convergence rate for this method is also very slow and is given by $\epsilon=O(N^{-1/2})$, where $N$ is the total number of points at which the deterministic model is solved. Another major issue is the lack of control of the distribution of points in the domain which causes unwanted clustering and scattering of points. For complicated deterministic models with high stochastic dimensions, the number of realizations required for a certain high level of accuracy may be unrealistic. Approaches like Latin Hypercube sampling (LHS) \cite{helton2003latin}, Importance Sampling \cite{rubinstein2011simulation,au2003important}, Quasi Monte Carlo Methods (using Halton sets, Sobol sets) \cite{niederreiter1992random} have been used successfully to achieve better convergence rates than the conventional Monte Carlo method. Artificial Neural Networks \cite{tu1996advantages}, In-situ Adaptive Tabulation (ISAT) \cite{pope1997computationally} and the Inverse Distance Weighted (IDW) \cite{shepard1968two} technique are some of the approaches which can be used in tandem with Monte Carlo sampling as postprocessing tools to approximate the surface and hence build a surrogate surface.\\
\indent Stochastic Galerkin Method \cite{le2010spectral} is a spectral approach which is a very popular tool for uncertainty propagation. It is a non-sampling approach where the unknown solution is projected onto the stochastic space spanned by a set of complete orthogonal polynomials after which the Galerkin projection is applied to minimize the error due to the gPC expansion and form a coupled set of deterministic equations. Wiener's original work on polynomial chaos \cite{wiener1938homogeneous} dealt with representation of a Gaussian random process using global Hermite polynomials. Initial work on the stochastic Galerkin method was done by Ghanem and Spanos \cite{ghanem2003stochastic} using the concept of polynomial chaos by Wiener and has been subsequently applied to various practical problems \cite{ghanem1998scales,ghanem1999stochastic,ghosh2008strain,doostan2009least,ghanem1999propagation}. Generalized polynomial chaos (gPC) expansion was developed by Xiu and Karniadakis \cite{xiu2002wiener} by including various other global polynomial-random variable combinations, with a few applications found in \cite{xiu2002modeling,xiu2003modeling,xiu2003new}. This method is known to have a very high convergence rate given the response surface is sufficiently smooth in all the stochastic dimensions. In the presence of discontinuities or highly localized variations in the response surface, this method may fail to converge due to the well-known Gibbs phenomenon. Remedies for this problem have been sought using multielement gPC \cite{kewlani2009multi,wan2005adaptive,wan2006multi}, piecewise polynomial basis \cite{babuska2004galerkin}, the wavelet basis \cite{le2004uncertainty} and basis enrichment of polynomial chaos expansions \cite{ghosh2008stochastic}.
All these methods involve solution of a coupled system of deterministic equations which may be non-trivial to solve when the original deterministic model is very complex in itself. This is the drawback of the intrusive nature of the method. A way to get around this issue is the usage of non-intrusive collocation approaches.\\
\indent The basic idea of non-intrusive collocation approaches is to strategically select points in the stochastic space. A surrogate response surface is then constructed based on these points to allow for cheap extraction of more samples. The goal is to achieve a specified level of accuracy with an optimally small number of sample evaluations. This method solves the deterministic problem at pre-selected collocation points in the random domain, determined by using either interpolation approaches or discrete projection approaches \cite{xiu2010numerical}. Some of the earlier works on this method \cite{babuvska2007stochastic,mathelin2003stochastic} used a tensor product of 1-D interpolation functions. This approach suffers from the so-called `curse of dimensionality' \cite{cools2002advances} as the number of points needed for full model evaluations increases exponentially with increase in the number of dimensions. Sparse grid \cite{bungartz2004sparse} approaches alleviate this problem to some extent as they significantly reduce the number of points in high dimensions while maintaining almost the same level of accuracy. Sparse grids are especially suitable for high dimensional problems involving numerical integration and interpolation. The interpolation approach approximates the stochastic space using multi-dimensional interpolation with the existing data such that the surrogate surface always passes through the pre-determined points \cite{smolyak1963quadrature,klimke2005algorithm,klimke2006uncertainty,klimke2007sparse}. More recent works introduce adaptivity into the sparse grid collocation interpolation approach, including dimension-adaptive sparse grid methods \cite{ganapathysubramanian2007sparse,nobile2008anisotropic}, Multi-Element(domain-adaptive) sparse grid interpolation \cite{foo2008multi,agarwal2009domain}, and adaptive sparse grid subset interpolation \cite{ma2009adaptive}. The adaptivity helps to efficiently characterizing any highly localized variations and discontinuities in the response surface. The discrete gPC projection approach, also known as the pseudospectral approach \cite{xiu2010numerical} is a discretized version of the exact generalized Polynomial Chaos(gPC) projection method, where a multi-dimensional numerical integration is performed with the existing data to approximate the stochastic solution. The surrogate surface here is not constrained to pass through the pre-determined points. This approach is non-intrusive and has fast convergence rates for smooth stochastic domains, but it is less conducive to tackling problems with discontinuous response in the stochastic space. A global approach based on Pad$\acute{e}$-Legendre approximation \cite{chantrasmi2009pade} has also been used to track down strong non-linearities or discontinuities in the response surface. It has also been shown \cite{elman2012stochastic,grigoriu2014response} that selection of input points by considering the probability structure of the input domain can lead to efficient sampling.\\
\indent The present work is based on the work done by Ma and Zabaras \cite{ma2009adaptive} on adaptive sparse grid subset interpolation. Similar to that work, the proposed approach uses linear basis functions for the adaptive sparse grid interpolation to capture any localized variations in the response. In addition, it aims to reduce the number of function evaluations by local 1-dimensional cubic spline interpolations \cite{kincaid2002numerical} in the smoother regions of the response domain. The smoothness is measured by successive derivative estimation along a straight line of points using finite differences of the output values in any of the input dimensions. Small changes (within a tolerance) in the derivative estimates will indicate sufficient smoothness for cubic spline interpolation along the straight line. This helps to achieve the same accuracy as in \cite{ma2009adaptive}, but decreases the number of function evaluations, especially when the response function is widely smooth. It is worth mentioning here that the derivative information is extracted approximately from the output values without any exact knowledge about the derivative of the output of interest.\\ 
\indent The rest of the manuscript is organized as follows: In section 2, the general mathematical model for any physical system with uncertainties is described. In section 3, the conventional stochastic collocation (CSC) method, the adaptive sparse grid collocation (ASGC) method and then the proposed efficient adaptive sparse grid collocation (E-ASGC) method are discussed in details. Section 4 deals with the various numerical examples to compare the performance of the proposed method with a few existing methods. Finally, the concluding remarks are given in section 5. \\ 

\section{Problem Definition}
Following notations in \cite{ma2009adaptive}, we represent the complete probability space by the triplet ($\Omega,\mathcal{F},\mathcal{P}$) where $\Omega$ corresponds to the sample space of outcomes, $\mathcal{F} \subset 2^{\Omega}$ is the sigma algebra of measurable events in $\Omega$, and $\mathcal{P} : \mathcal{F} \rightarrow [0,1] $ is the probability measure.
Let $I(\omega)=\{I_1, I_2, I_3 ,….., I_d\}$ be the multidimensional vector of random input parameters in a problem of interest, where $ I : \Omega \rightarrow \in \mathbb{R}^d $ 
\begin{equation}
Z(\omega)=f(I(\omega)), \ \ \forall \  \omega \in \Omega
\end{equation}
The goal is then to find out how the vector valued output $Z(\omega)$ varies with respect to each of the random vector components $I_i(\omega), i \in [1,2,$..$,d]$.\\
\section{Stochastic Collocation Interpolation Method}
\subsection{Conventional Sparse Grid Interpolation}For a function $f:[a,b]\rightarrow \mathbb{R}$, the one-dimensional interpolation formula is given by:
\begin{equation}
U^k(f(x))=\sum_{x^k \in X^k}a_{x^k}(x)f(x^k)=\sum_{j=1}^{m_k}a_{x^k_j}(x)f(x^k_j),\\
\end{equation}
where $x \in [a,b],$\ $X^k=\{x^k | x^k \in [a,b]\},$\ $a_{x^k_j}(x) \in [0,1] \subset \mathbb{R}^1, \ a_{x^k_j}(x^{k}_i)=\delta_{ij},$ \  $\{i,j\} \in [1,2,...m_k],$ and  $m_k$ = number of points in the set $X^{k}.$ For multi-dimensional interpolation, the one-dimensional case can be upgraded to obtain a tensor product formulae:
\begin{equation}
(U^{k_1} \otimes ...... \otimes  U^{k_d} )(f(\mathbf{x}))=\sum_{j_1=1}^{m_1}....\sum_{j_d=1}^{m_d}(a_{x^{k_1}_{j_1}}(\mathbf{x})\otimes...... \otimes a_{x^{k_d}_{j_d}}(\mathbf{x}))f(x^{k_1}_{j_1},....,x^{k_d}_{j_d})\\
\end{equation}
where $d$ is the total number of dimensions and $\mathbf{x}=\{x_1,x_2,...,x_d \} \in \mathbb{R}^d$\\
\indent The major drawback of this tensor product formula is that the total number of points required are $(m_1)(m_2)(m_3)......(m_d)$ which rises exponentially with increase in dimensions, leading to the curse of dimensionality. The sparse grid approach that is used in the current work mitigates this issue to quite an extent by sampling significantly fewer points which are subsets of the tensor grid structure. Though the accuracy of the algorithm is not totally dimension-independent, it gets weakened down to a logarithmic dependence.\\
\indent Using similar definitions as in \cite{klimke2005algorithm}, we define $U^0=0$ and the incremental interpolant by:
\begin{equation}
\triangle^k(f(x))=U^k(f(x))-U^{k-1}(f(x)), \ \ \forall k \ge 1 
\end{equation}
where,
\begin{equation}
U^k(f(x))=\sum_{x^k \in X^k}a_{x^k}(x)f(x^k)\\
\end{equation}
and
\begin{equation}
U^{k-1}(f(x))=U^k(U^{k-1}(f(x)))
\end{equation}
By using the above three equations, we thus get,
\begin{align}
\triangle^k(f(x)) &  =\sum_{x^k \in X^k}a_{x^k}(x)f(x^k)-\sum_{x^k \in X^k}a_{x^k}(x)g(x^k) \nonumber \\
& =\sum_{x^k \in X^k}a_{x^k}(x)(f(x^k)-g(x^k))
\end{align}
\ \ \ \ \ \ \ \ \ \ \ \ \ \ where $g=U^{k-1}(f(x))$\\
Now, 
\begin{equation} 
(f(x^k)-g(x^k))=0,  \ \ \ \forall x^k \in X^{k-1}
\end{equation}
Thus,
\begin{equation}
\triangle^k(f(x))= \sum_{x^k \in X_{\triangle}^{k-1}}a_{x^k}(x)(f(x^k)-g(x^k))
\end{equation}
where $X_{\triangle}^k=X^k \setminus X^{k-1} $ denotes the points in set $X^k$ but not in $X^{k-1}$. Because of the nested property of uniform grids, $X^{k-1} \subset X^k$,
the number of elements (points) in $X^{k}_{\triangle}=m_k-m_{k-1}=m^k_{\triangle}$\\
\indent Rewriting Eq. (9), we get,
\begin{equation}
\triangle^k(f(x))=\sum_{j=1}^{m^k_{\triangle}}a_{x^k_j}(x)(f(x^k_j)-g(x^k_j))
\end{equation}
Using the property $\triangle^k(f(x))=U^k(f(x))-U^{k-1}(f(x))$,we can write
\begin{equation}
U^{k}(f(x))=\sum_{i=1}^k\triangle^i(f(x))\\
\end{equation}
In the case of a tensor grid, the multivariate interpolant expression is a tensor product extension of Eq. (11) and is given by,
\begin{equation}
(U^{k_1} \otimes U^{k_2} \otimes ....... \otimes  U^{k_d} )(f(\mathbf{x}))=\sum_{i_1=1}^{k_1}.....\sum_{i_d=1}^{k_d}(\triangle^{i_1}\otimes.......\otimes\triangle^{i_d})(f(\mathbf{x}))
\end{equation}
On the other hand, Smolyak sparse grids \cite{smolyak1963quadrature} use a much smaller subset of the tensor grid. The sparse grid interpolant only considers points satisfying $|\textbf{i}| \le q $ and is defined by,
\begin{equation}
A_{q,d}=\sum_{|\textbf{i}|\le q} (\triangle^{i_1} \otimes ......\otimes \triangle^{i_d})(f(\mathbf{x}))=A_{q-1,d}(f(\mathbf{x}))+\triangle A_{q,d}(f(\mathbf{x}))
\end{equation}
where $|\textbf{i}|=i_1+i_2+.....i_d,$ for $\textbf{i}=(i_1,i_2,....i_d) \in \mathbb{N}^d$ and $q=k+d-1$ for $k_1=k_2=\dots=k_d=k \ge 1$.
\begin{align}
\triangle A_{q,d}(f(\mathbf{x}))&= \sum_{|\textbf{i}|=q} (\triangle^{i_1} \otimes ......\otimes \triangle^{i_d})(f(\mathbf{x}))\\
A_{q-1,d}(f(\mathbf{x}))&= \sum_{|\textbf{i}| \le q-1} (\triangle^{i_1} \otimes ......\otimes \triangle^{i_d})(f(\mathbf{x}))
\end{align}
Here $i_s$ $\left\{ \forall s=1,2,.....,d \right\}$  is called the depth of interpolation in the $s$-th dimension whereas q denotes the global depth of interpolation and $ A_{d-1,d}=0$.
Thus, putting $q=d$, we get $A_{d,d}(f)=\triangle A_{d,d}(f(\mathbf{x}))$ and for $q=d+1, $ $A_{d+1,d}(f(\mathbf{x}))=A_{d,d}(f(\mathbf{x}))+\triangle A_{d+1,d}(f(\mathbf{x}))$ and so on. Thus it is seen that there is a hierarchy when it comes to the interpolant at different levels, and the interpolant at a given level contributes to the estimation of the next higher level interpolant.\\
\indent From Eqs. (10) and (14), we get,
\begin{align}
\triangle A_{q,d}(f(\mathbf{x}))&= \sum_{|\textbf{i}|=q} (\triangle^{i_1} \otimes ......\otimes \triangle^{i_d})(f(\mathbf{x})) \nonumber \\
&= \sum_{|\textbf{i}|=q} \sum_{\mathbf{j}} (a^{i_1}_{j_1}(\mathbf{x})\otimes.......\otimes a^{i_d}_{j_d}(\mathbf{x}))(f(x^{i_1}_{j_1},.....,x^{i_d}_{j_d}) \nonumber \\
&\qquad {} -
g_1(x^{i_1}_{j_1},.....,x^{i_d}_{j_d}))\nonumber 
\end{align}
where $\mathbf{j}=\{(j_1,j_2,...,j_d): j_s=1,...,m^{i_s}_{\triangle};\  s=1,...,d\}$, and $g_1$ is defined as
\begin{align}
g_1 &=(U^{k_1-1} \otimes ....... \otimes  U^{k_d-1} )(f(\mathbf{x}))\nonumber\\
&= \sum_{|\mathbf{i}|\le q-1} (\triangle^{i_1} \otimes ......\otimes \triangle^{i_d})(f(\mathbf{x})) \nonumber\\
&=A_{q-1,d}(f(\mathbf{x})) \nonumber
\end{align}
Now,
\begin{align}
\triangle A_{q,d}(f(\mathbf{x}))&= \sum_{|\textbf{i}|=q} \sum_{\mathbf{j}} (a^{i_1}_{j_1}(\mathbf{x})\otimes.......\otimes a^{i_d}_{j_d}(\mathbf{x}))(f(x^{i_1}_{j_1},.....,x^{i_d}_{j_d}) \nonumber \\
&\qquad {} -
g_1(x^{i_1}_{j_1},.....,x^{i_d}_{j_d})) \nonumber \\
&= \sum_{|\textbf{i}|=q} \sum_{\mathbf{j}}a^{\mathbf{i}}_{\mathbf{j}}(\mathbf{x})w^{\mathbf{i}}_{\mathbf{j}}(\mathbf{y})
\end{align}
where $a^{\mathbf{i}}_{\mathbf{j}}$ is the $d$-dimensional basis function, $w^{\mathbf{i}}_{\mathbf{j}}$ is the hierarchical surplus and $\mathbf{y}=\{x^{i_1}_{j_1},.....,x^{i_d}_{j_d}\}$.\\
Thus, once the surrogate model has been identified for any given level $q$, the function value at any point can be calculated as:
\begin{equation}\label{eq_surrogate}
u(\mathbf{x},\mathbf{y})= \sum_{|\mathbf{i}| \le q} \sum_{\mathbf{j}}a^{\mathbf{i}}_{\mathbf{j}}(\mathbf{x})w^{\mathbf{i}}_{\mathbf{j}}(\mathbf{y})
\end{equation}
The mean of the solution can be analytically estimated \cite{ma2009adaptive} as:
\begin{equation}
\mathbb{E}[u(\mathbf{x},\mathbf{y})]=\bar{u}(\mathbf{y})=\int_{\Gamma}\sum_{|\mathbf{i}| \le q} \sum_{\mathbf{j}}w^{\mathbf{i}}_{\mathbf{j}}(\mathbf{y}) a^{\mathbf{i}}_{\mathbf{j}}(\mathbf{x})\rho(\mathbf{x})d \mathbf{x}
\end{equation}
Since $\mathbf{x}$ is an uniform random space, the probability density function $\rho(\mathbf{x})=1$ for the domain $\Gamma=[0,1]^d.$
Substituting this value of $\rho(\mathbf{x})$, rearranging the integral and assuming the random variables are independent of each other,
\begin{equation}
\bar{u}(\mathbf{y})=\sum_{|\mathbf{i}| \le q} \sum_{\mathbf{j}}w^{\mathbf{i}}_{\mathbf{j}}(\mathbf{y}) l^{\mathbf{i}}_{\mathbf{j}}
\end{equation}
where
\begin{equation}
l^{\mathbf{i}}_{\mathbf{j}} = \int_{\Gamma} a^{\mathbf{i}}_{\mathbf{j}}(\mathbf{x})d \mathbf{x}=\prod_{k=1}^{d} \int_{0}^{1} a^{i_k}_{j_k}(x)dx \end{equation}
and  
\[
    \int_{0}^{1} a^{i_k}_{j_k}(x)dx= 
\begin{cases}
    1,					& \text{if } i_k = 1\\
    \frac{1}{4},      & \text{if } i_k =2\\
    2^{1-i_k},  	    & \text{otherwise}\
\end{cases}
\]\\
To get the variance, we proceed with the square of the function value in Eq. (\ref{eq_surrogate}) which can be approximated as:
\begin{equation}
u^2(\mathbf{x},\mathbf{y})= \sum_{|\mathbf{i}| \le q} \sum_{\mathbf{j}}a^{\mathbf{i}}_{\mathbf{j}}(\mathbf{x})v^{\mathbf{i}}_{\mathbf{j}}(\mathbf{y})
\end{equation}
where $v^{\mathbf{i}}_{\mathbf{j}}$ is the hierarchical surplus corresponding to the square of the output.\\
Then the expectation of the square of the random solution can be estimated as:
\begin{equation}
\mathbb{E}[u^2(\mathbf{x},\mathbf{y})]=\bar{u^2}(\mathbf{y})=\sum_{|\mathbf{i}| \le q} \sum_{\mathbf{j}}v^{\mathbf{i}}_{\mathbf{j}}(\mathbf{y}) \int_{\Gamma} a^{\mathbf{i}}_{\mathbf{j}}(\mathbf{x})d \mathbf{x}
\end{equation}
Thus the variance of the solution is given by:
\begin{equation}
Var[u(\mathbf{x},\mathbf{y})]=\sigma^2_u(\mathbf{y})=\bar{u^2}(\mathbf{y})-\bar{u}(\mathbf{y})^2
\end{equation}
\subsection{Adaptive Sparse Grid Interpolation}In conventional sparse grid interpolation methods, the error check is such that if any hierarchical surplus at the current level of interpolation exceeds the tolerance, all points in the next higher level must be evaluated. The algorithm ignores the fact that there may be smooth regions which do not require subsequent refinements. In adaptive sparse grid interpolation, unnecessary higher level samples are avoided by performing selective refinements. This method makes use of the tree-like data structure of 1-D equidistant sparse grid points. A schematic of the adaptive procedure using the tree-like structure of the grid points is given in Figure 1. With the exception of the point addition at level 2, two points are added in the neighborhood of each point at the previous level. In \cite{ma2009adaptive}, a point at the current level has been referred to as a `father' while points added around it at the next level are referred to as `sons'. For a $d$-dimensional random domain, there will be $2d$ sons added for each father if it is not a level-1 point in any of the dimensions. Therefore, if there are $m$ points at the current level at which the hierarchical surplus exceeds the tolerance, then at most $2dm$ points are added at the next level. There may be duplication of next-level points, requiring that a check be performed to avoid redundant sampling. This approach leads to slow convergence in the presence of localized variations or discontinuities, requiring computation up to a very high interpolation level. On the other hand, regardless of the nature of the response surface, the mean and variance converge quickly to a desired level of accuracy because of the sharp drop in the integral weights with rise in interpolation level. So it is reasonable to limit the algorithm to a maximum interpolation level which acts as another termination criterion.\\
\begin{figure}
\centering
\includegraphics[width=1\linewidth, height=4cm]{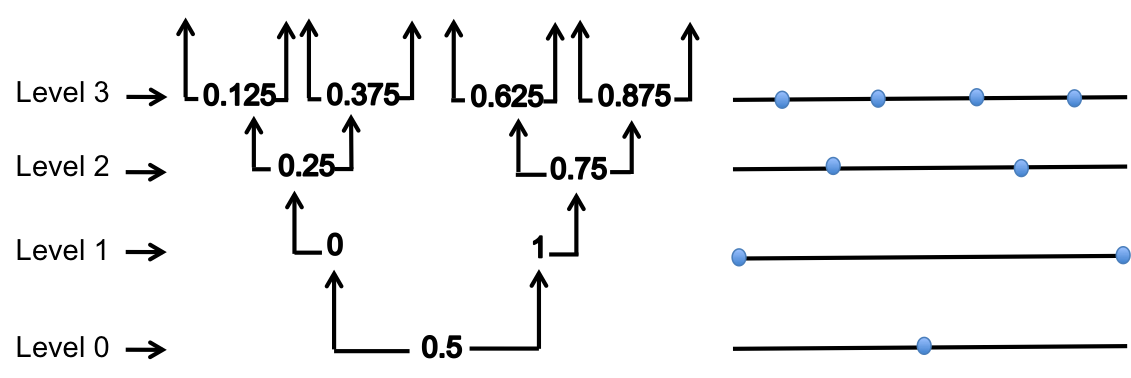} 
\caption{Tree-like structure of 1-D grid points}
\label{fig:subim0131}
\end{figure}
\indent An example of the adaptive procedure in 1-D is given in Figure 2. A 1-D $C^1$ discontinuous function shown in Figure 2(a) is used to demonstrate the adaptivity of the method. Figure 2(b) shows the conventional approach of adding the sparse grid points up to level 5. The absence of any adaptivity leads to a total of 33 point evaluations. In contrast, the adaptive collocation strategy shown in Figure 2(c) allows for local refinement of points around the $C^1$ discontinuous region and results in only 17 function evaluations up to level 5. For example, at level 3, out of the four points, only two of the points at $x=0.125$ and $x=0.375$ are `fathers' to points in the next higher level. This implies that the hierarchical surpluses at those two points are above the tolerance, while the hierarchical surpluses at $x=0.625$ and $x=0.875$ are below the tolerance.

\begin{figure}
\centering
\begin{subfigure}[b]{0.5\textwidth}
\centering
\includegraphics[width=1.2\linewidth, height=4cm]{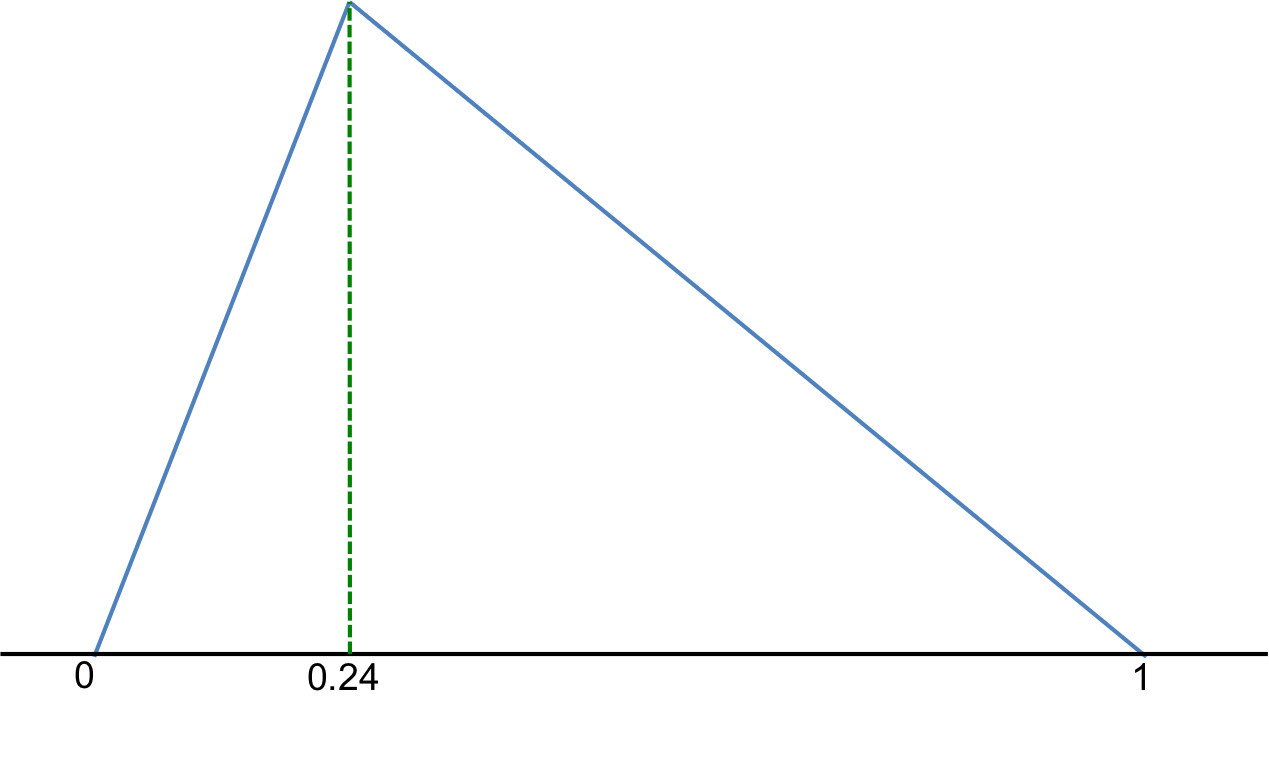} 
\caption{A 1-D function with $C^1$ discontinuity}
\label{fig:subim151}
\end{subfigure}
\\
\hfill 
\begin{subfigure}[b]{1\textwidth}
\includegraphics[width=1\linewidth, height=4cm]{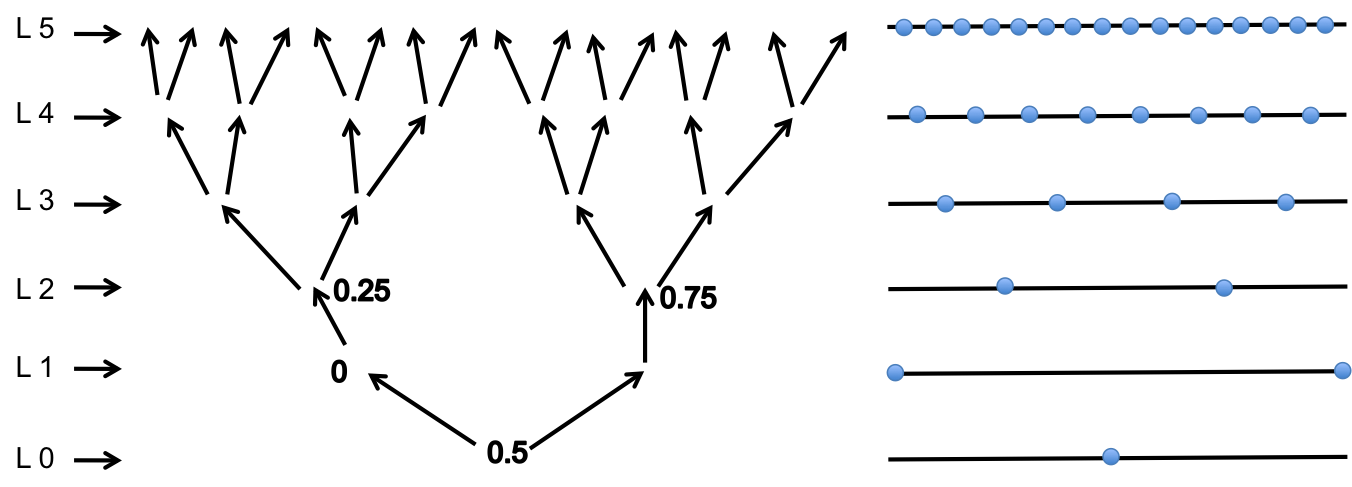} 
\caption{Conventional 1-D sparse grid points}
\label{fig:subim131}
\end{subfigure}
\\
\hfill 
\begin{subfigure}[b]{1\textwidth}
\includegraphics[width=1\linewidth, height=4cm]{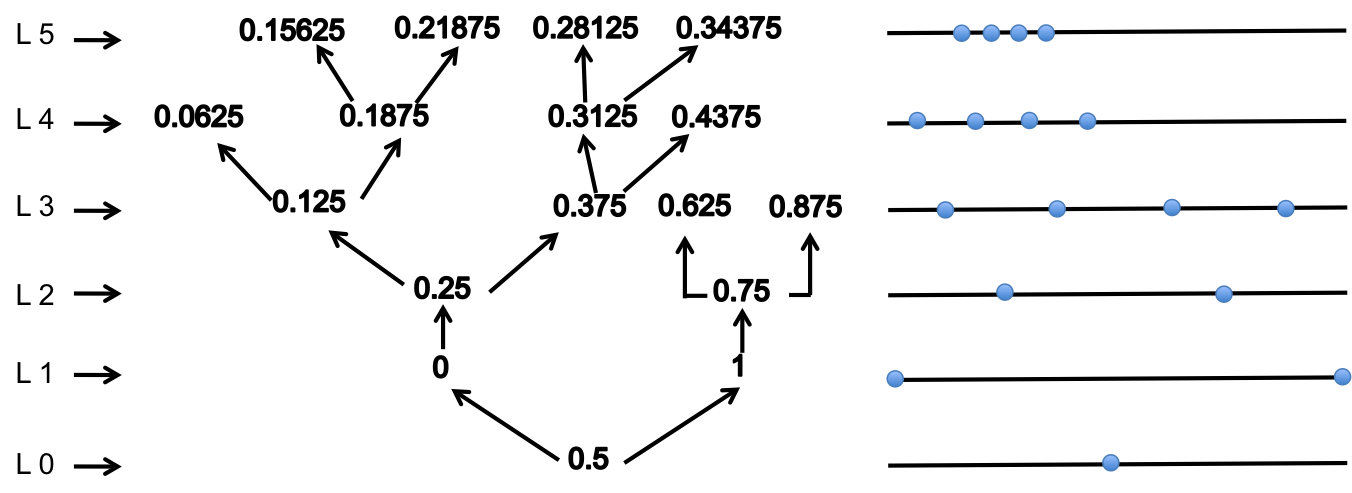} 
\caption{Adaptivity on 1-D sparse grid points}
\label{fig:subim131}
\end{subfigure}
\caption{Comparison between the addition of the conventional and the adaptive sparse grid points in 1-D}
\label{fig:image3}
\end{figure}
\indent In a high dimensional case with highly localized variations along some dimensions, there may be a significant number of other dimensions along which the response function is smooth without any sharp variations. Thus using a piece-wise linear function leads to slow convergence of the surpluses in the smooth regions. A significant number of full model evaluations can be avoided by handling these smoother regions efficiently. The proposed efficient  collocation method is based on this very idea which will be discussed next.\\
\subsection{Efficient adaptive sparse grid collocation through derivative estimation} One aspect to improve in the adaptive sparse grid subset collocation algorithm \cite{ma2009adaptive} is that in allocating more points along discontinuities and important dimensions, there can be a significant number of points also added to smoother regions in the domain. These points that are assigned to the smoother regions unnecessarily increase the computational cost. Thus one way to improve efficiency is to avoid brute force evaluations in the smoother regions as much as possible. The proposed method aims at achieving this by approximating the smoother regions with cubic splines \cite{kincaid2002numerical}. Therefore when proceeding with the adaptive algorithm \cite{ma2009adaptive}, if a sparse grid node is generated within any 1-D approximated smooth region, the function value at that point can be approximated by cubic spline interpolation. Cubic splines, being third order polynomials, can achieve sufficiently fast convergence and are robust because of their piece-wise nature. Higher order polynomials were not used to avoid over-fitting.\\

\begin{figure}
\centering
\includegraphics[width=0.9\linewidth, height=9cm]{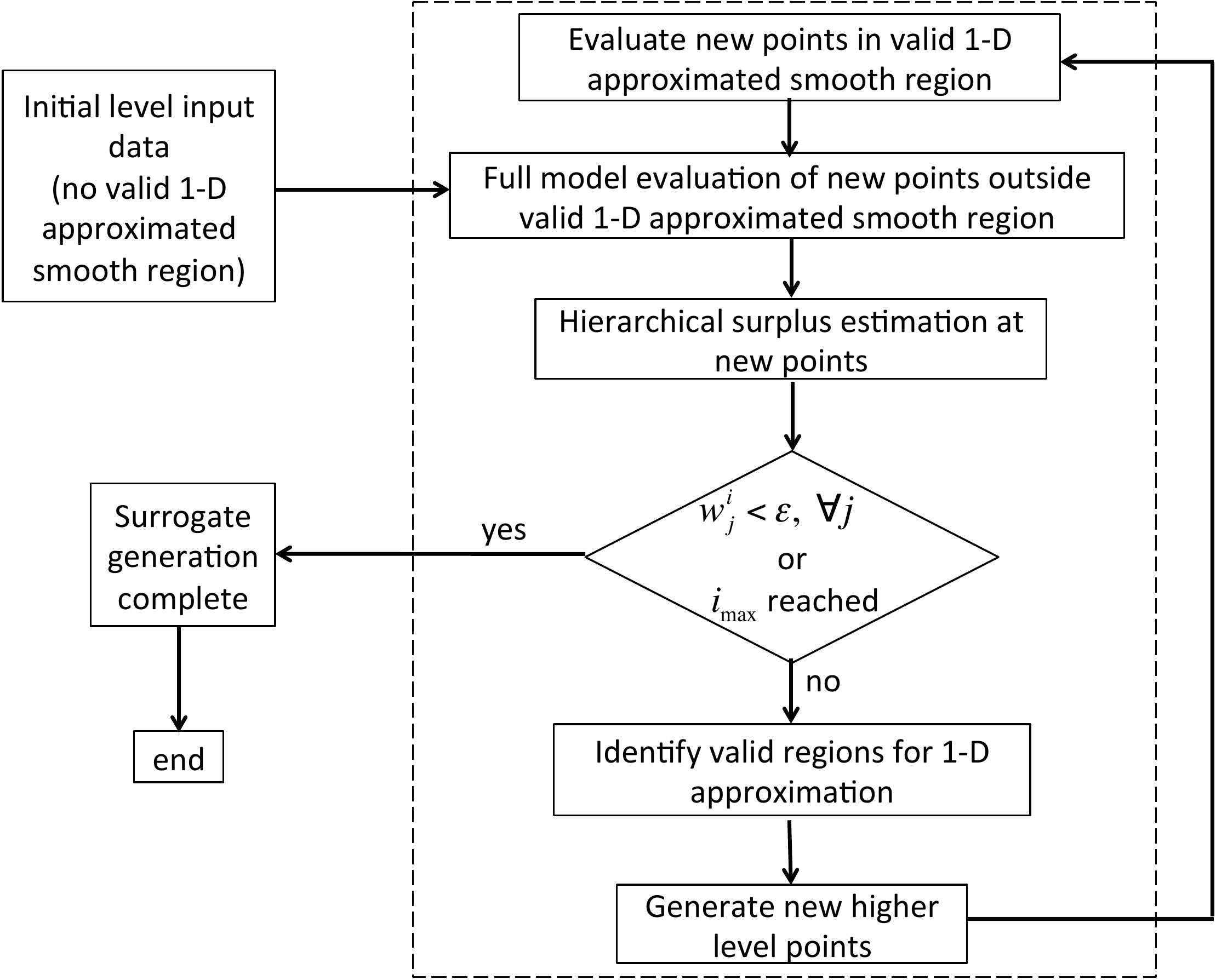} 
\caption{{Flow chart for the efficient adaptive sparse grid collocation algorithm}}
\label{fig:subim031}
\end{figure}

\indent The mechanism works as follows: Let P be the total number of unique sparse grid points and $D_i \ [ i \in {1,2,....d} ]$ be a certain dimension in the $d$-dimensional stochastic space. A projection is now taken on the plane orthogonal to the $D_i$ dimension to form a $D_{d-1}$ dimensional non-unique data point set. This set is non-unique because of the overlap of multiple points due to the projection. The critical control parameter is the minimum number of points in a single straight line along  the $D_i$ dimension needed to approximate the linear region. The non-unique set is used to evaluate the frequency of points along every straight line in the $D_i$ dimension. Once the criterion is met, successive finite difference based derivative calculation is done to crudely detect any discontinuity along the line. If the change in the derivative varies gradually enough throughout, the data along that path is stored. Otherwise the line can be split into 1-D sub-intervals and data in those individual sub-intervals are stored. The data storage include the input points, the function value at the points, the midpoint of the linear interval and a parameter $L$ quantifying the extent of the approximation along a straight line. $L$ is defined as half of the total length of any 1-D approximated smooth region. In future calculations, if we come across a point lying in this 1-D region where the function value needs to be known, then instead of performing an expensive brute force evaluation there, we retrieve the associated data and simply approximate the function value by the cubic spline interpolant.\\
The algorithm is outlined below:\\
\underline{Initialization}\\
\indent (1) Set the maximum interpolation level $i_{max}$\\
\indent (2) Set the tolerance parameter $\epsilon$\\
\indent (3) Set the dimension $d$ of the problem\\
\indent (4) Initialize the level to $i=1$\\
\underline{Conventional sparse grid iterations for initial few levels $i_1$}\\
\indent While $i \le i_1$,\\
\indent (5) Perform full model evaluation at level $i$ points\\
\indent (6) Calculate also the hierarchical surpluses $w_j^i$ at each of the points\\
\indent (7) Set $i=i+1$\\
\indent (8) Form $2d$ points at $(i+1)$-th level \\
\indent (9) Check for duplication of previous points\\ 
\indent (10) Go to step (5)\\
\underline{Main adaptive loop based on derivative estimation}\\
\indent While $i_1 < i \le i_{max}$,\\
\indent (11) Check if any of the points lie in the regions identified as smooth enough for 1-D interpolation by searching the stored database.\\
\indent (12) For each point lying in any of these regions, perform a cheap interpolation to get the function value.\\
\indent (13) For the remaining points, perform full model evaluations.\\
\indent (14) Calculate the hierarchical surpluses $w_j^i$ at each of the points in steps (12) and (13)\\
\indent (15) Check if $|w_j^i| \ge \epsilon$ for each point. If no points satisfy this criterion, go to step (22).\\
\indent (16) For each point satisfying criterion in step (15), form $2d$ points at the $(i+1)$-th level\\
\indent (17) Set $i=i+1$\\
\indent (18) Check for duplicity of points\\ 
\indent (19) For each dimension, project all points on the hyperplane normal to the dimension, to get $M_p$ number of $(d-1)$-dimensional non-unique projected points. Count the number of unique points $N$ and the number of co-located points at each of the unique points $m_i$, such that $\sum_i^N m_i =M_p$, where $m_i \ge 1$\\
\indent (20) For each point $\in M_p$ with $m_i$ greater than a minimum number $M_{min}$,\\
\indent  \indent	a) Calculate successive derivative along that straight line using finite differences of the output values.\\
\indent  \indent	b) If the change in the gradients are below a tolerance $\phi$, the interval is considered smooth, and it is stored in the database for future retrieval. The database for an interval consists of the input data, the output data, the midpoint value of the interval, and half the length of the interval, represented by the quadruplet $(I,O,\bar{I},L)$.\\
\indent (21) Go to step (11)\\
\indent (22) End generation of surrogate model.\\
\indent The generated surrogate model can be used as the basis for generating moments (e.g., mean and variance) and/or distribution of the response. Also, the output response at any arbitrary query point can be extracted. The flow chart for the entire algorithm is shown in Figure 3.\\
\\
\textbf{Remark}: It is very important to note here that poor choices of parameters $M_{min}$ and $\phi$ may render the proposed method very inefficient in certain situations. If $M_{min}$ is not large enough, then the finite difference derivatives could be inaccurate, which could lead to inaccurate use of the cubic spline interpolation, particularly if the tolerance $\phi$ is too high. If the cubic spline interpolated value at a sparse grid point is erroneous as a result of these poorly chosen parameters, then one of two undesirable scenarios could occur. First, the hierarchical surplus error at that point could appear to be larger than the tolerance when in reality it is not, which would direct the algorithm to add new and unnecessary sampling points around that point. Second, the hierarchical surplus error at the point could appear to be smaller than the tolerance when in reality it is not. This would direct the algorithm to add sampling points further away from this point, potentially missing local variations in the response function near this point and reducing the accuracy of the results. If these scenarios are encountered quite often in a particular surrogate modeling procedure, then the performance of the proposed method will be worse than the ASGC method both in terms of the error and the number of function evaluations at a particular sparse grid level.\\
\subsection{Convergence and accuracy of the proposed efficient adaptive method}In the proposed adaptive method based on 1-D derivative estimation, the accuracy depends heavily on the minimum number of points $M_{min}$ chosen for the storage and retrieval of points for cubic spline interpolation.
\begin{equation}
\|u^q_{asgc}-u^q_{e-asgc}\|_\infty  \le \epsilon_1 N_a
\end{equation} 
where $\epsilon_1\le \frac{5}{384} ||u^{q (4)}_{asgc}||_\infty h^4$  is the interpolation error \cite{hall1976optimal} which decreases with increase in $M_{min}$, $h$ is the maximum knot spacing and $N_a$  are the points where full model evaluations are performed in the ASGC \cite{ma2009adaptive} sparse grid but are interpolated in the E-ASGC sparse grid. In the interpolation error expression, $(.)^{(4)}$ denotes the fourth order derivative. \\
Now, the interpolation error \cite{ma2009adaptive} between the adaptive sparse grid and the conventional sparse grid is given by:
\begin{equation}
\|u^q_{csc}-u^q_{asgc}\|_\infty \le \epsilon_2N_b
\end{equation}
where $N_b$ are the points in the conventional sparse grid but missing in the adaptive sparse grid subset due to the hierarchical surplus based adaptivity.\\
The interpolation error of the conventional sparse grid \cite{klimke2005algorithm} is given by:
\begin{equation}
\| u-u^q_{csc}\|_\infty = O(N^{-2}_t |log_2N_t|^{3(d-1)})
\end{equation}
where $N_t$= total number of interpolation points at interpolation depth $q$ in $d$- dimensional stochastic space.\\
Thus the approximate bound of the total error using the proposed method is:
\begin{align}
\|u-u^q_{e-asgc}\|_\infty&=\| u-u^q_{csc}+u^q_{csc}-u^q_{asgc}+u^q_{asgc}-u^q_{e-asgc} \|_\infty \nonumber\\
& \le \| u-u^q_{csc}\|_\infty +\| u^q_{csc}-u^q_{asgc}\|_\infty +\| u^q_{asgc}-u^q_{e-asgc} \|_\infty
\end{align}
It is worth mentioning that $M_{min}$ should be chosen such that $\epsilon_1 < \epsilon_2$.
\section{Numerical Examples}In this section, results are shown for explicit functions as well as implicit functional variations in different dimensions. All examples compare the adaptive sparse grid interpolation \cite{ma2009adaptive} and the proposed efficient adaptive sparse grid collocation method (E-ASGC). The first example studies 2-dimensional analytic functions. The second example considers a family of 5-dimensional analytic functions. In the third example, a spatial one-dimensional elliptic problem with high-dimensional stochasticity is used to compare the performance of the different methods. The last example deals with an indeterminate truss structure with variable cross-sectional areas of its members.\\
\subsection{Function with $C^1$ discontinuity}We consider the function in $[0,1]^2$ as mentioned in \cite{ma2009adaptive}
\begin{equation}
f(x,y)=\frac{1}{|0.3-x^2-y^2|+0.1}
\end{equation}
The exact function is plotted in Figure 4(a) and it is seen that there is a line discontinuity which is not along any of the two dimensions. It is also observed that away from the discontinuities, the function is quite smooth which is suited for higher order interpolation. The proposed derivative based approach aims to utilize this feature of the function.
\begin{figure}
\centering
\begin{subfigure}[b]{0.4\textwidth}
\centering
\includegraphics[width=1.3\linewidth, height=5cm]{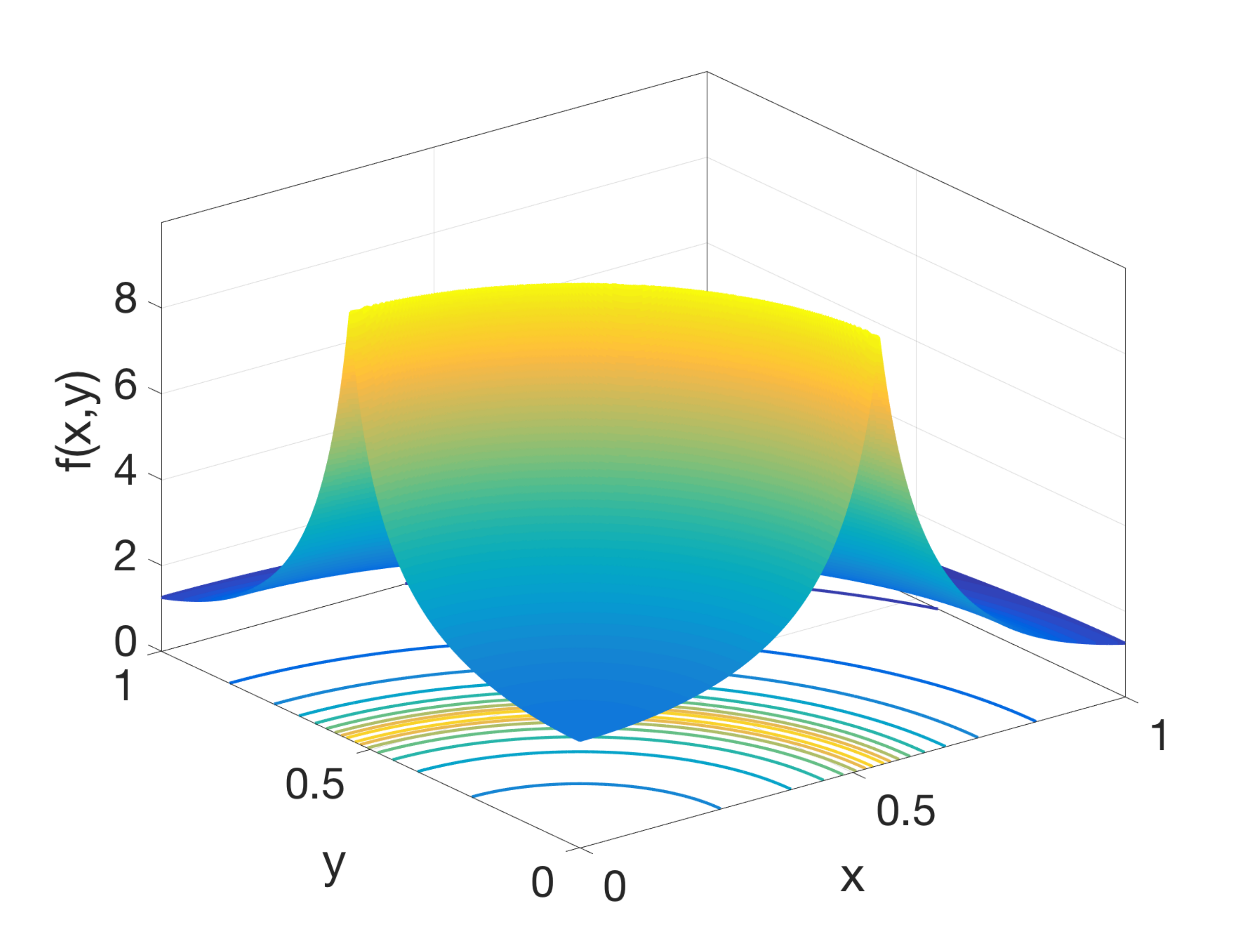} 
\caption{}
\label{fig:subim231}
\end{subfigure}
\hfill
\begin{subfigure}[b]{0.4\textwidth}
\includegraphics[width=1.3\linewidth, height=5cm]{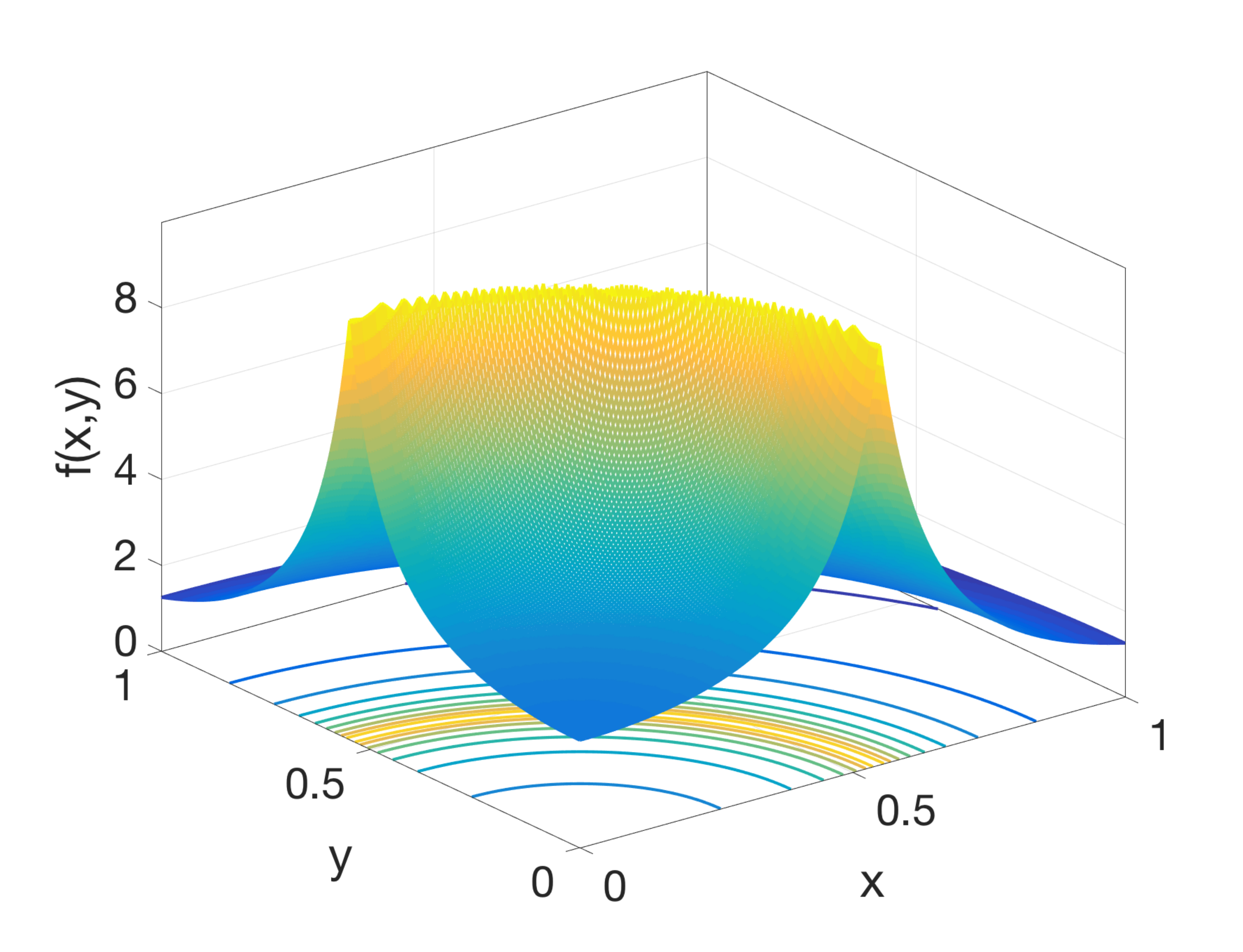}
\caption{}
\label{fig:subim241}
\end{subfigure}
\hfill 
\begin{subfigure}[b]{0.4\textwidth}
\centering
\includegraphics[width=1.3\linewidth, height=5cm]{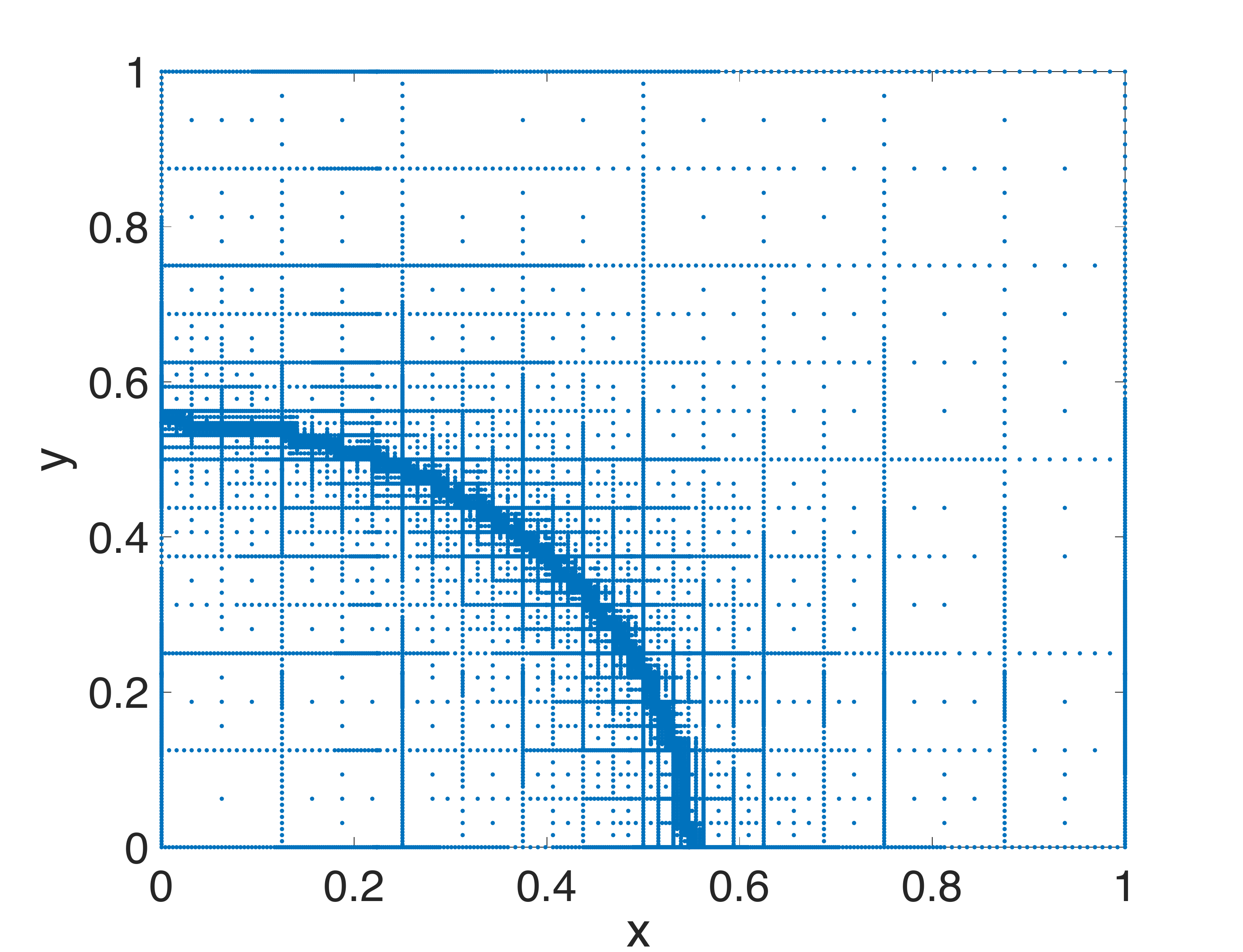}
\caption{}
\label{fig:subim241}
\end{subfigure}
\hfill 
\begin{subfigure}[b]{0.4\textwidth}
\centering
\includegraphics[width=1.3\linewidth, height=5cm]{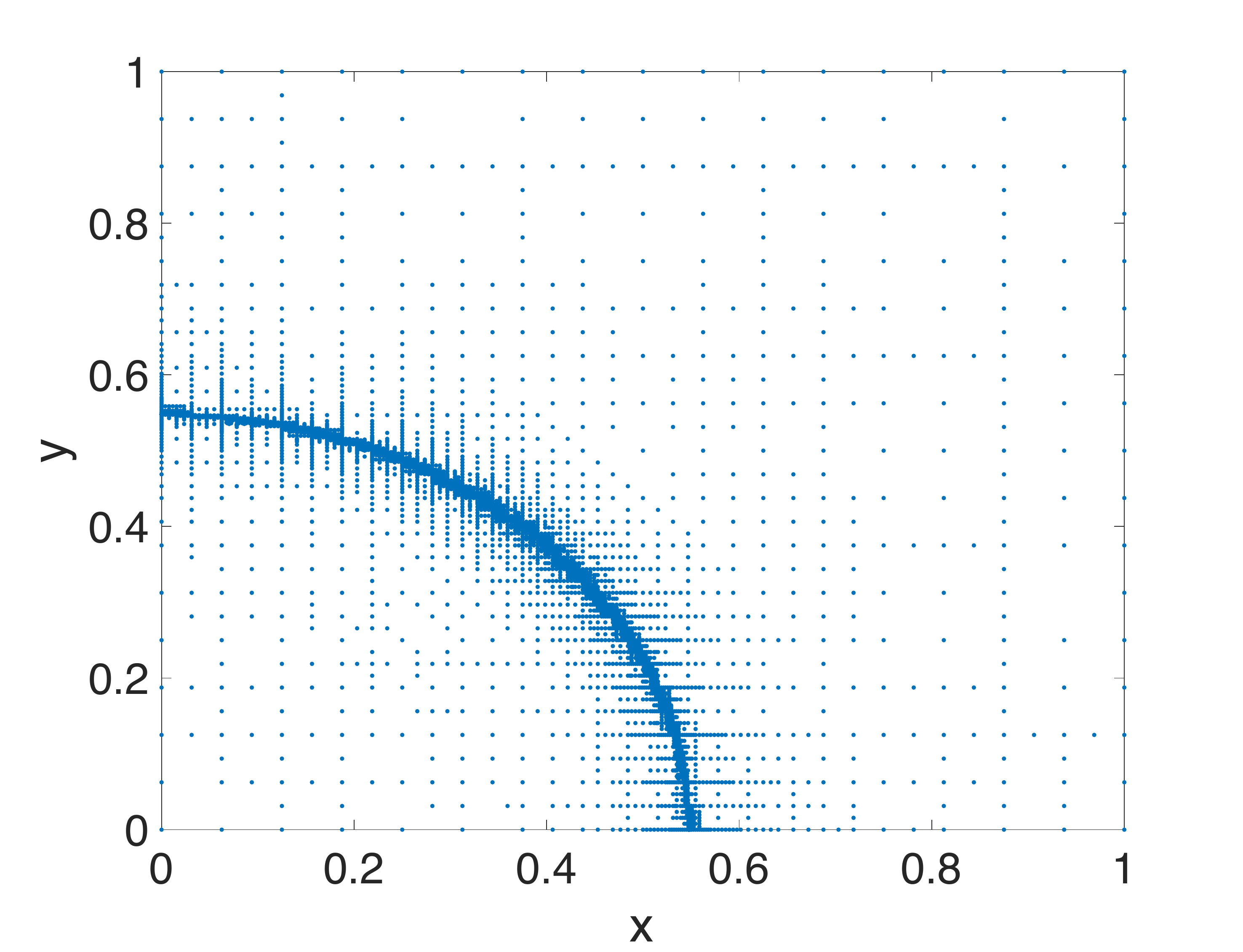}
\caption{}
\label{fig:subim241}
\end{subfigure}
\hfill 
\begin{subfigure}[b]{0.4\textwidth}
\centering
\includegraphics[width=1.3\linewidth, height=5cm]{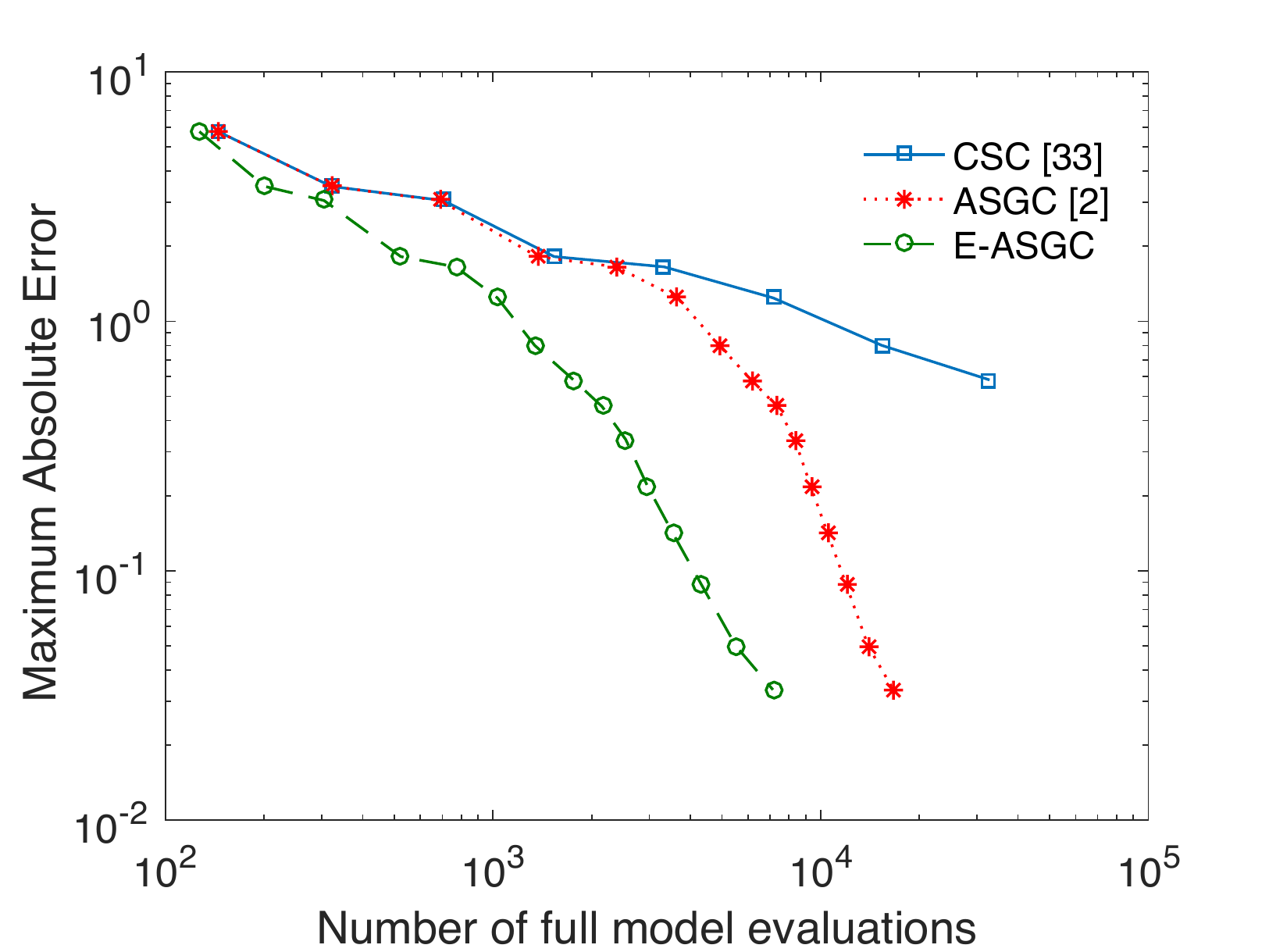}
\caption{}
\label{fig:subim251}
\end{subfigure}
\hfill 
\begin{subfigure}[b]{0.4\textwidth}
\centering
\includegraphics[width=1.3\linewidth, height=5cm]{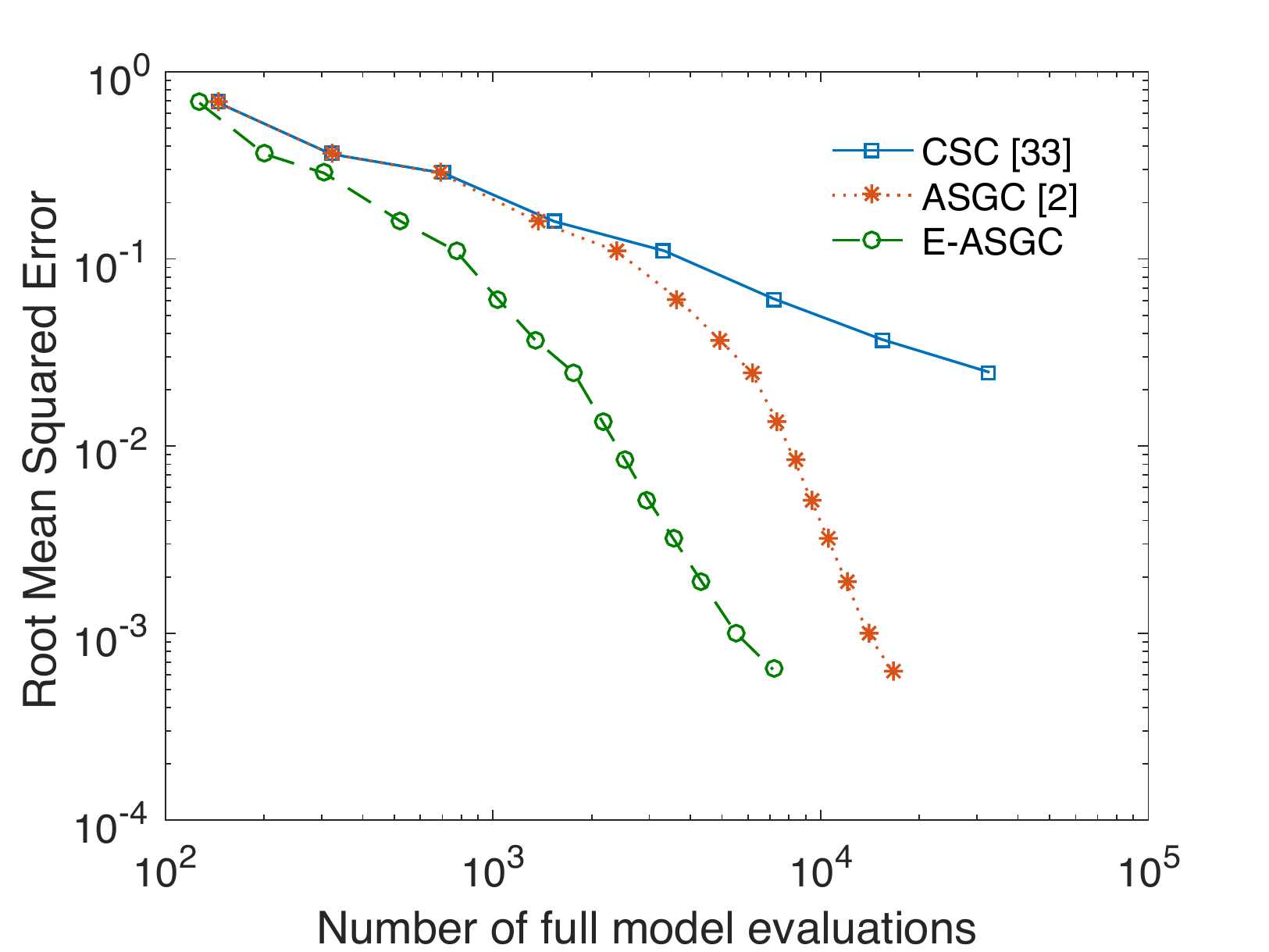}
\caption{}
\label{fig:subim251}
\end{subfigure}
\caption{(a) Exact 2D line singularity function; (b) Approximate function from E-ASGC; (c) ASGC sampling points (16,659); (d) E-ASGC sampling points (7,149) (e) Maximum Absolute error convergence plot comparing conventional stochastic collocation (CSC), adaptive sparse grid subset collocation (ASGC) and efficient adaptive sparse grid collocation (E-ASGC); (f) Root Mean Squared Error convergence plot}
\label{fig:image4}
\end{figure}
The approximate function obtained from the E-ASGC algorithm at the interpolation depth 19 is shown in Figure 4(b). The sampling points for the existing ASGC algorithm are shown in Figure 4(c) while those of the E-ASGC method are shown in Figure 4(d). It can seen from the plots that the E-ASGC approach effectively approximates the smooth regions of the input domain, with significantly fewer full model evaluations in those regions. Figures 4(e) and 4(f) shows the error convergence plots for the conventional sparse grid collocation (CSC), ASGC and E-ASGC methods. It is seen that the E-ASGC method clearly outperforms both the ASGC and the CSC methods. As a comparison, for a maximum absolute error of 0.0334, the E-ASGC method requires 7,149 function evaluations while the ASGC method requires 16,659 function evaluations. Thus the E-ASGC method reduces sampling by more than a factor of 2 relative to the ASGC method. The CSC method has the worst performance with a total of 32,769 function evaluations required for a maximum absolute error of 0.5824.\\
\subsection{5 dimensional functions}The 5-dimensional family of functions are taken from Genz \cite{genz1987package}. They were primarily used to assess the efficiency of numerical integration algorithms. The functions defined on $x \in [0,1]^5$ are described \cite{klimke2005algorithm} as follows:\\
Oscillatory function:
\begin{equation}
f_1(\mathbf{x})=\mathrm{cos}(2 \pi w_1+\sum_{i=1}^{5}c_i x_i)
\end{equation}
where $w_1$ and $c_i : i=1,2,....5$ are constants.\\
Corner Peak Function:
\begin{equation}
f_2(\mathbf{x})=(1+\sum_{i=1}^{5}c_i x_i)^{-6}
\end{equation}
where $c_i : i=1,2,....5$ are constants.\\
Discontinuous function:\\
\[
    f_3(\mathbf{x})= 
\begin{cases}
    0,& \text{if } x_1 \geq w_1 \text{ or } x_2 \ge w_2,\\
    \exp(\sum_{i=1}^{5}c_i x_i),              & \text{otherwise}
\end{cases}
\]
where $w_1,w_2$ and $c_i : i=1,2,....5$ are constants.\\
Convergence plots of the maximum absolute error for the E-ASGC method is compared with that of the conventional sparse grid method and the ASGC method for all the three functions, shown in Figure 5. These results show that the E-ASGC and ASGC are at least as efficient as the CSC method in all the cases and they perform significantly better for the corner peak and discontinuous functions. E-ASGC reduces the required samples from ASGC somewhat in all the three cases.\\ 
\begin{figure}
\centering
\begin{subfigure}[b]{0.4\textwidth}
\centering
\includegraphics[width=1.5\linewidth, height=6cm]{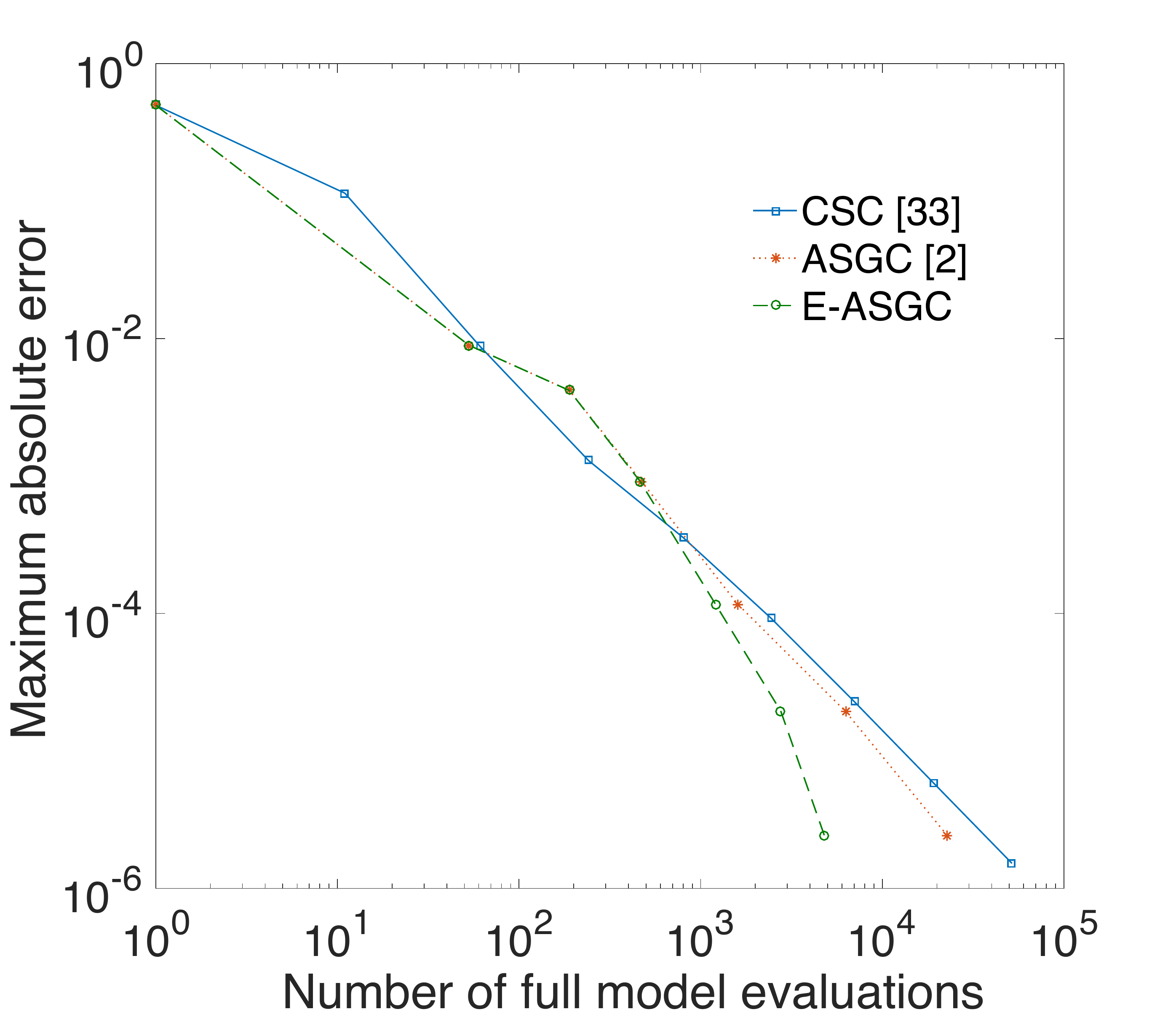} 
\caption{Oscillatory function}
\label{fig:subim151}
\end{subfigure}
\hfill 
\begin{subfigure}[b]{0.4\textwidth}
\centering
\includegraphics[width=1.5\linewidth, height=6cm]{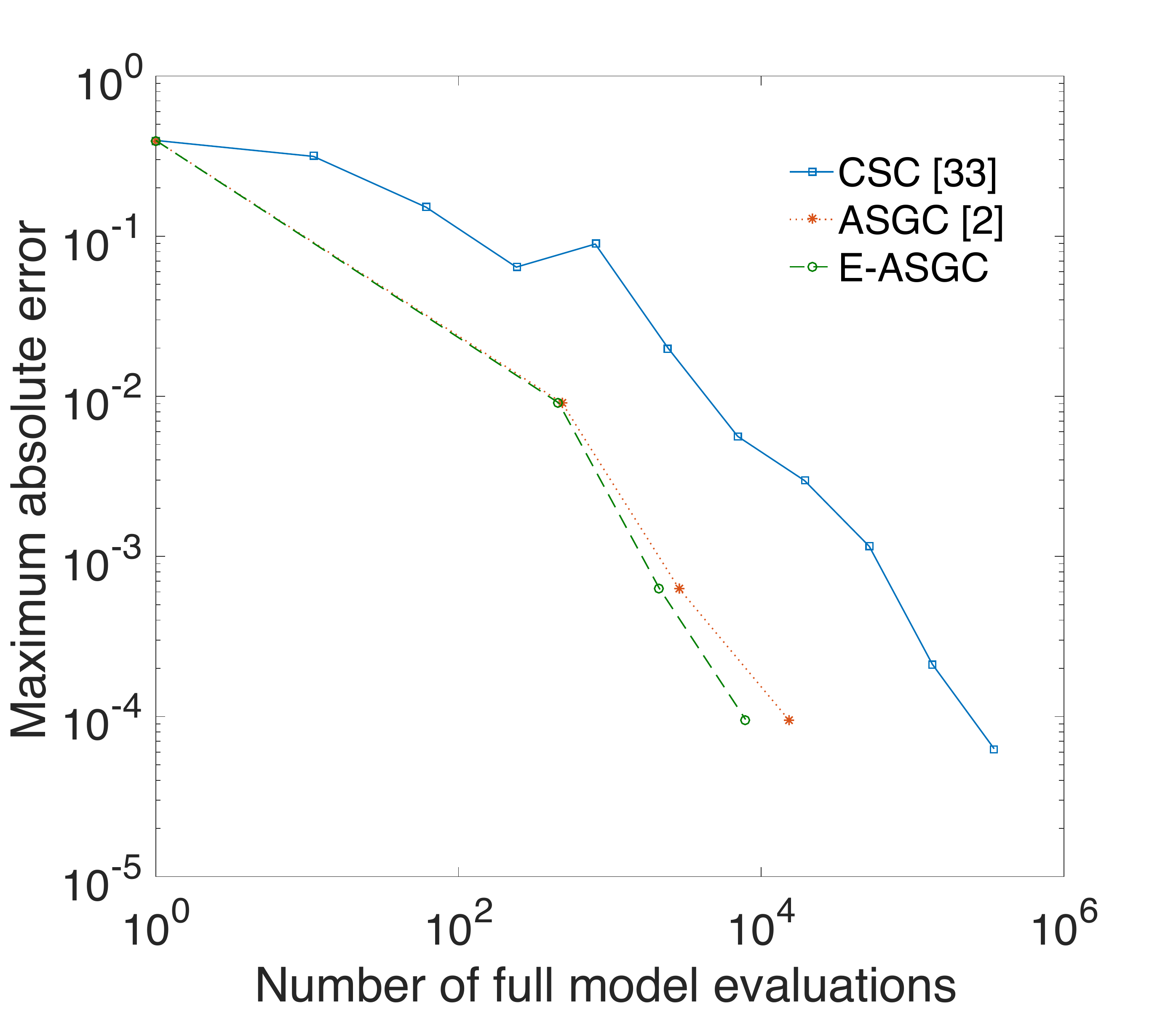} 
\caption{Corner Peak function}
\label{fig:subim131}
\end{subfigure}
\hfill
\begin{subfigure}[b]{0.4\textwidth}
\centering
\includegraphics[width=1.5\linewidth, height=6cm]{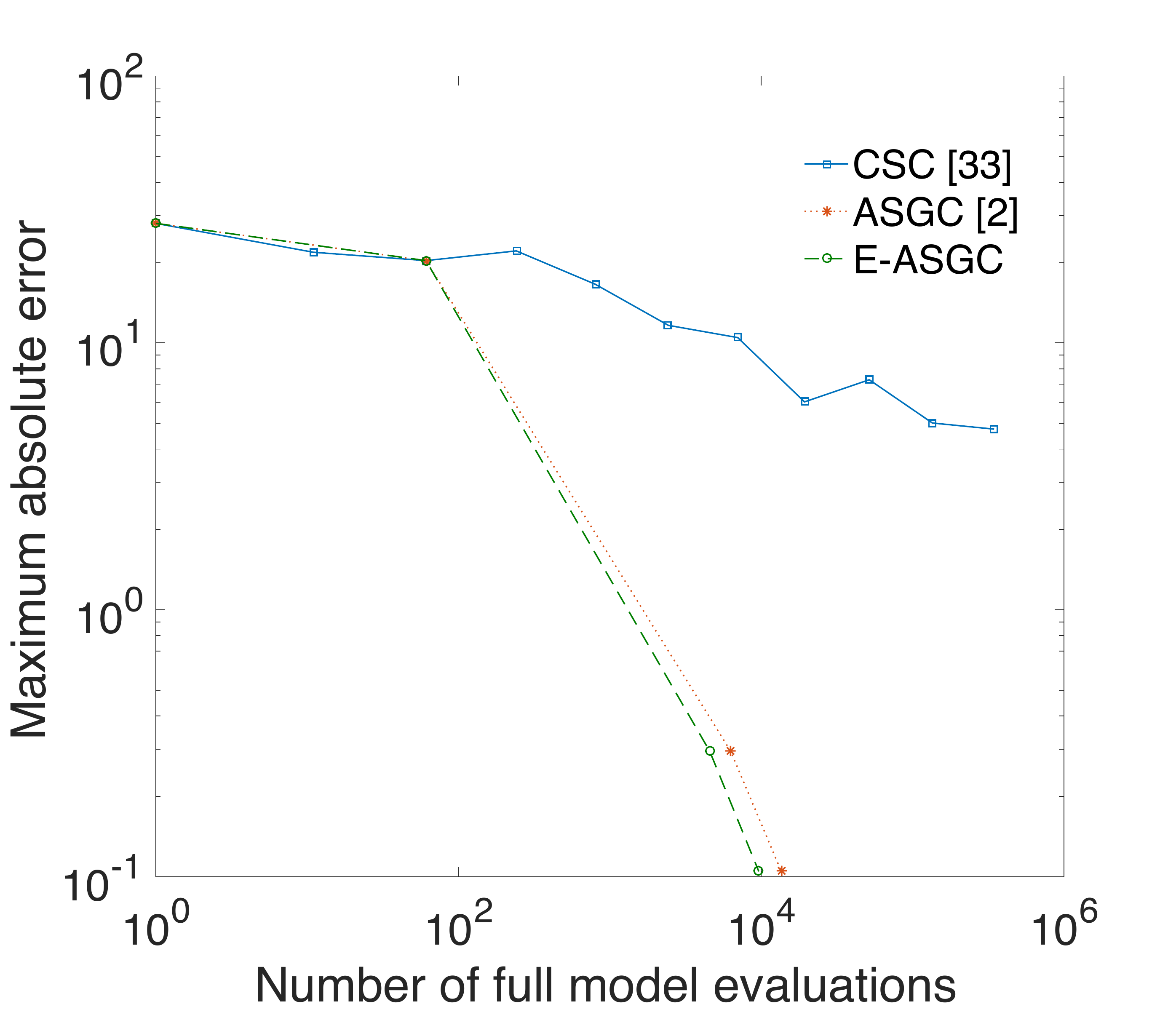}
\caption{Discontinuous function}
\label{fig:subim141}
\end{subfigure}
\caption{Error estimation for different 5 dimensional analytic functions}
\label{fig:image3}
\end{figure}

\subsection{Spatial 1-D Poisson Equation}This problem deals with a random process, making it an infinite-dimensional stochastic problem. The model problem is given by:
\begin{equation}
-\triangledown (\kappa(x,\omega) \triangledown u (x,\omega))=f(x)
\end{equation}
\begin{equation}
u(0,\omega)=u(1,\omega)=0
\end{equation}
\begin{equation}
f(x)=2x
\end{equation}
The diffusion coefficient $\kappa(x,\omega)$ is represented by a random process \cite{nobile2008anisotropic} and is approximated as a finite-dimensional random quantity by:\begin{equation}
\log(\kappa(x,\omega)-0.5) \approx 1+Y_1(\omega)(\frac{\sqrt{\pi}L}{2})^{1/2}+\sum_{i=2}^{N} \xi_n \phi_n(x) Y_n(\omega) 
\end{equation}
where,
\begin{equation}
\xi_n=(\sqrt{\pi}L)^{1/2} \exp(\frac{-( \lfloor {\frac{n}{2}}\rfloor  \pi L)^2}{8}), \text{ if } n>1
\end{equation}
and\\
\[
    \phi_n(x):= 
\begin{cases}
    \mathrm{sin}(\frac{\lfloor {\frac{n}{2}}\rfloor  \pi x}{L_p}),			& \text{if n is even},\\
    \mathrm{cos}(\frac{\lfloor {\frac{n}{2}}\rfloor  \pi x}{L_p}),              & \text{if n is odd}
\end{cases}
\]
where $Y_n(\omega)$ $ \{n=1,2,3......N\}$ are independent random variables which are uniformly distributed in $[0,1]$, $L_p=max{[1,2L_c]}$, $L=\frac{L_c}{L_p}$, where $L_c$ is the correlation length of the random process.\\ 
\begin{figure}[htbp]
\centering
\begin{subfigure}[b]{0.4\textwidth}
\centering
\includegraphics[width=1.45\linewidth, height=6cm]{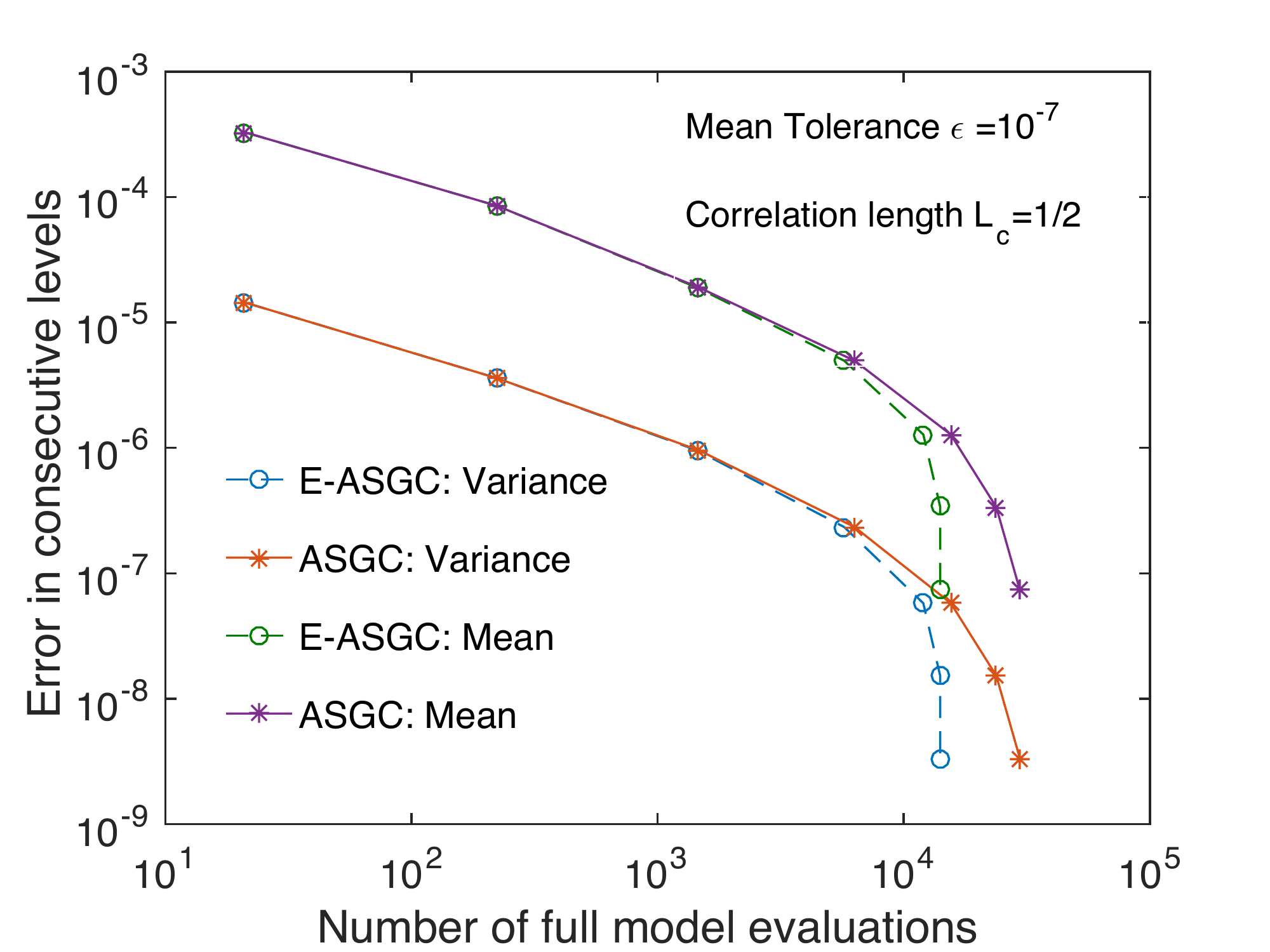} 
\caption{10 dimensions}
\label{fig:subim331}
\end{subfigure}
\hfill
\begin{subfigure}[b]{0.4\textwidth}
\centering
\includegraphics[width=1.45\linewidth, height=6cm]{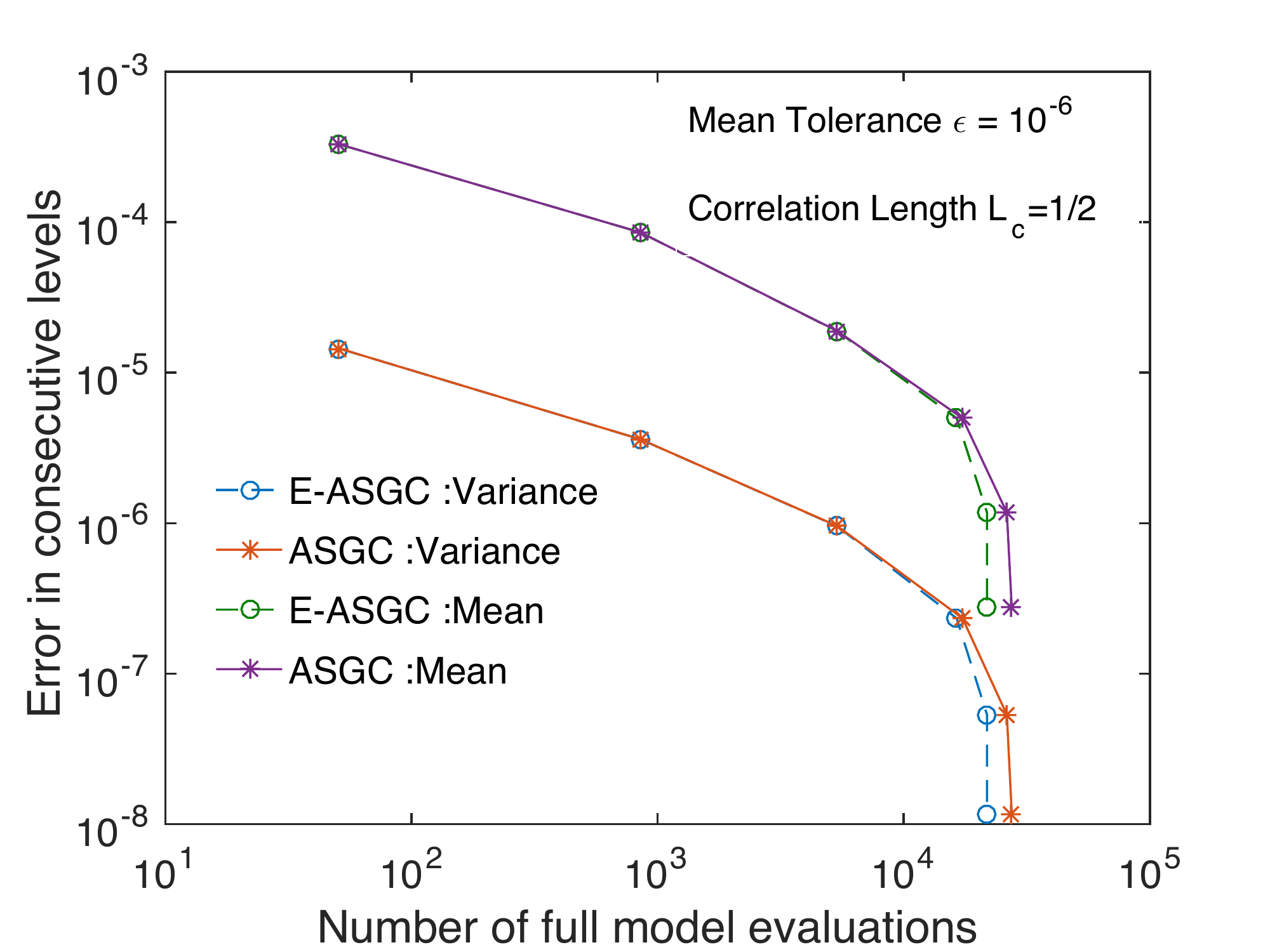}
\caption{25 dimensions}
\label{fig:subim341}
\end{subfigure}
\hfill 
\begin{subfigure}[b]{0.4\textwidth}
\centering
\includegraphics[width=1.45\linewidth, height=6cm]{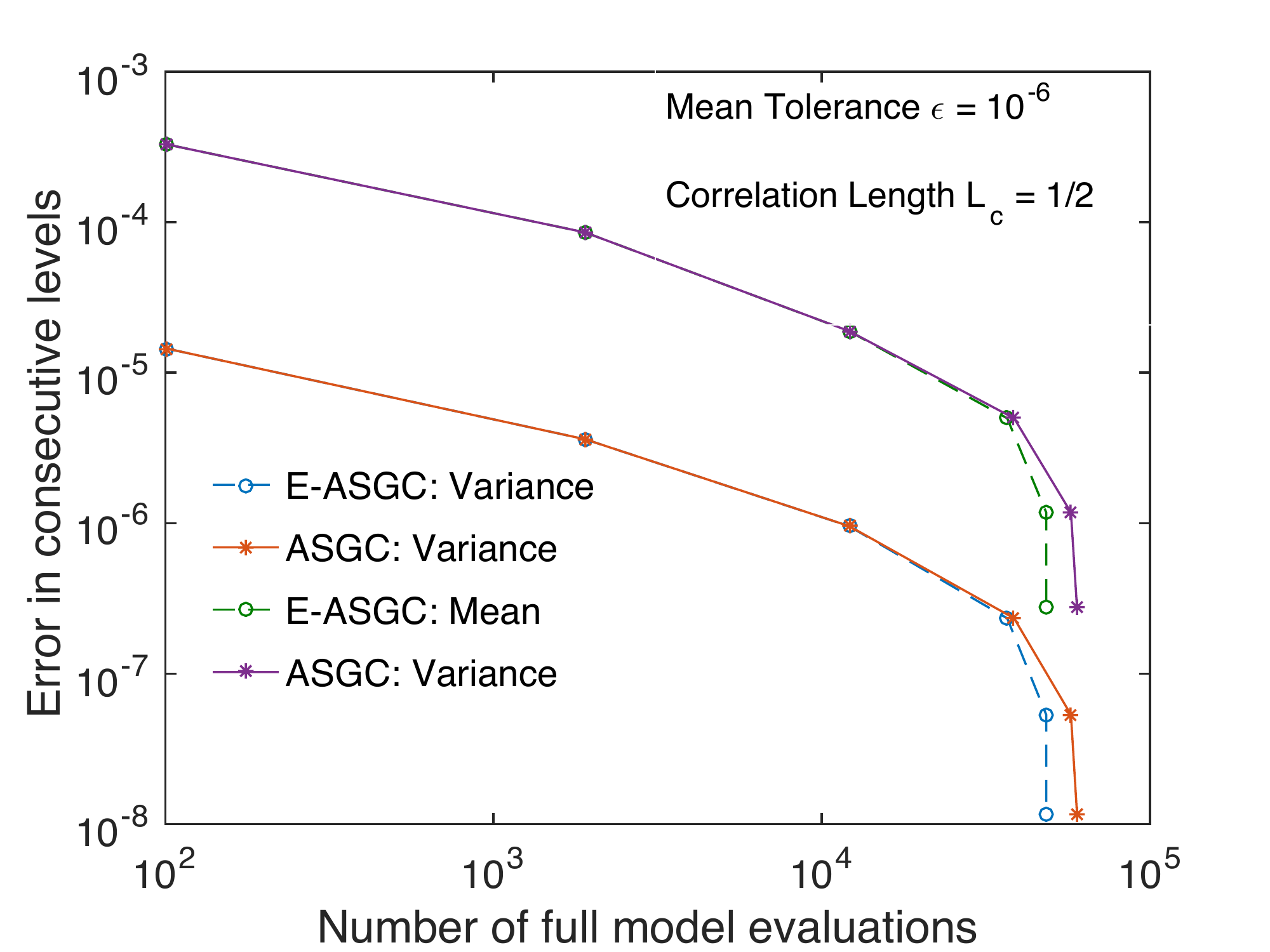} 
\caption{50 dimensions}
\label{fig:subim351}
\end{subfigure}
\hfill 
\begin{subfigure}[b]{0.4\textwidth}
\centering
\includegraphics[width=1.45\linewidth, height=6cm]{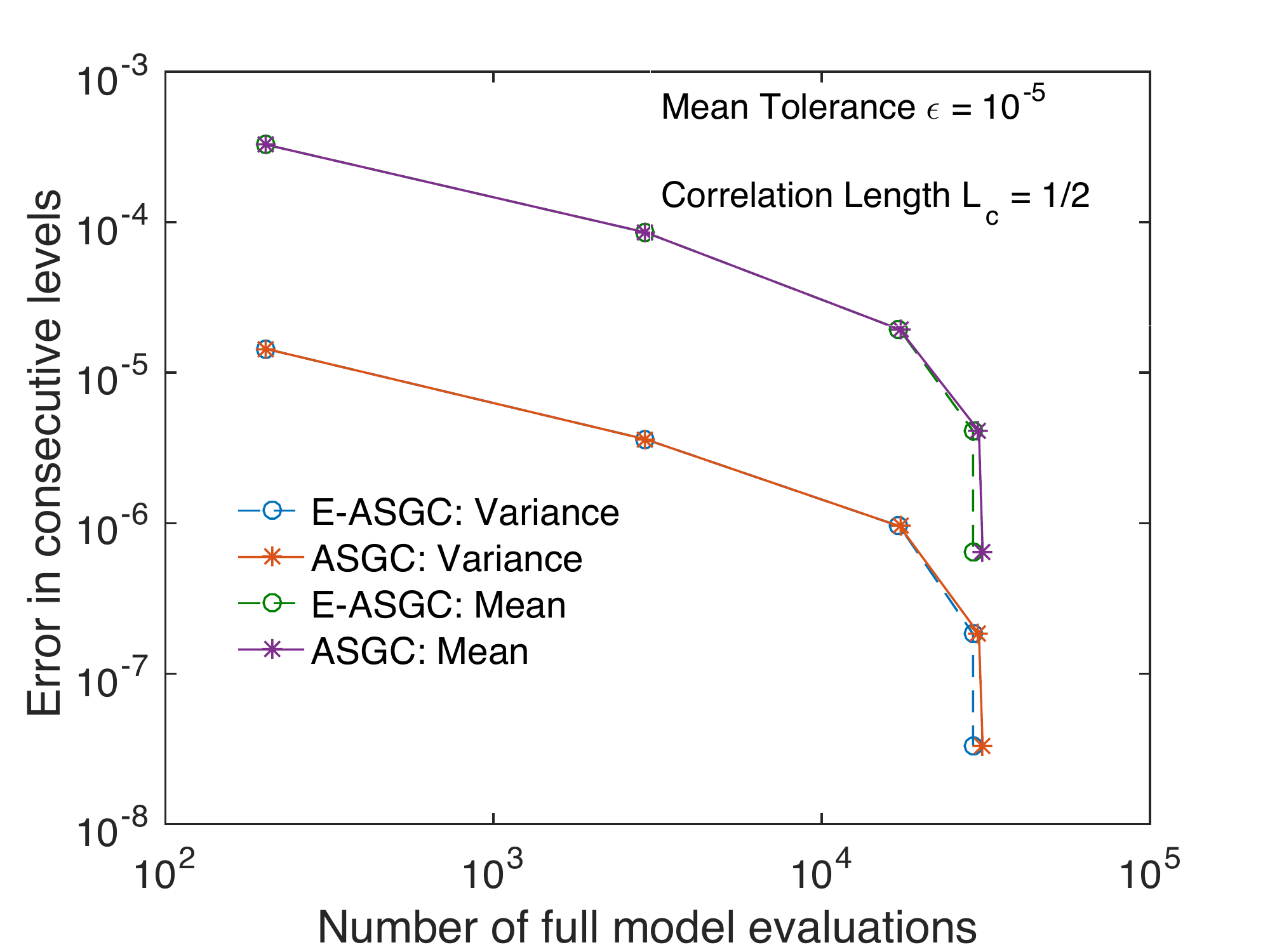} 
\caption{100 dimensions}
\label{fig:subim351}
\end{subfigure}

\caption{Error in Mean and Variance for different dimensions in the spatial 1-D stochastic elliptic problem}
\label{fig:image4}
\end{figure}
\indent A high value of correlation length implies that the eigenvalues decay fast and the first few dimensions are of more importance than the rest of the dimensions. To make this method self-sufficient for error estimation, the error for the mean between consecutive levels is given by $E(A_{q,N}(u_N))-E(A_{q,N+1}(u_N))$ which progressively converges to the target precision. Similarly the error in variance is given by $Var(A_{q,N}(u_N))-Var(A_{q,N+1}(u_N))$. The convergence in the mean and variance for this problem using the efficient adaptive sparse grid collocation (E-ASGC) method is shown and compared with the ASGC method in Figure 6. The results show that for relatively low accuracy, the ASGC and E-ASGC have the same efficiency in all the four cases. With increase in the accuracy, it is seen that the E-ASGC performs more efficiently than the ASGC approach although the efficiency seems to decrease with increase in the dimensionality of the problem from $N=10$ to $N=100$.\\

\subsection{Truss Problem}
We consider the 2D truss structure shown in Figure 7(a). It is assumed that all of the elements have the same modulus, E. The lengths of the vertical elements, horizontal elements and the diagonal elements are $\sqrt{3}L$, $L$ and $2L$ respectively. The truss is statically indeterminate to the first degree and is analyzed using the stiffness method. The uncertainty lies in the cross-sectional area of the truss members. The variations in the cross-sectional area of the members are such that for certain combinations of the areas, the force in member 5 exceeds its critical buckling load. It is then considered to have failed and carries no load. In that scenario, the truss effectively converts to that shown in Figure 7(b). The output of interest of this problem is the force in member 4 under the given loads, as a function of the member cross-sectional areas. If member 5 fails, it leads to a discontinuous increase in the force in member 4. The indeterminate truss structure and the determinate structure without the diagonal cross-brace member 5 are shown in Figure 7.
Node C in the truss structure is subjected to a vertical load $\sqrt{3}P$ acting downwards and a horizontal load $P$ to the left.

\begin{figure}
\centering
\begin{subfigure}[b]{0.4\textwidth}
\centering
\includegraphics[width=1.4\linewidth, height=6cm]{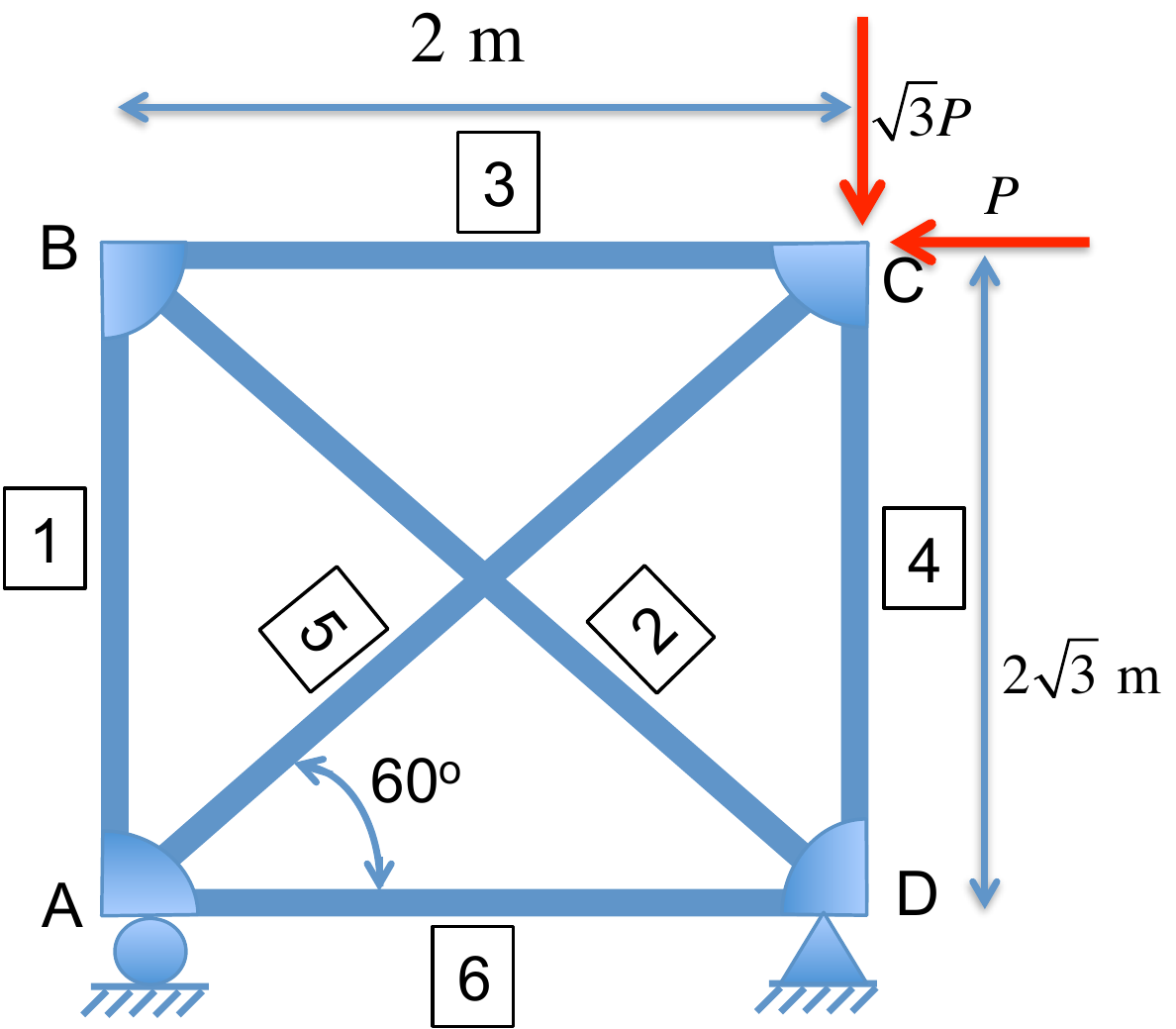} 
\caption{Statically indeterminate}
\label{fig:subim331}
\end{subfigure}
\hfill
\begin{subfigure}[b]{0.4\textwidth}
\centering
\includegraphics[width=1.4\linewidth, height=6cm]{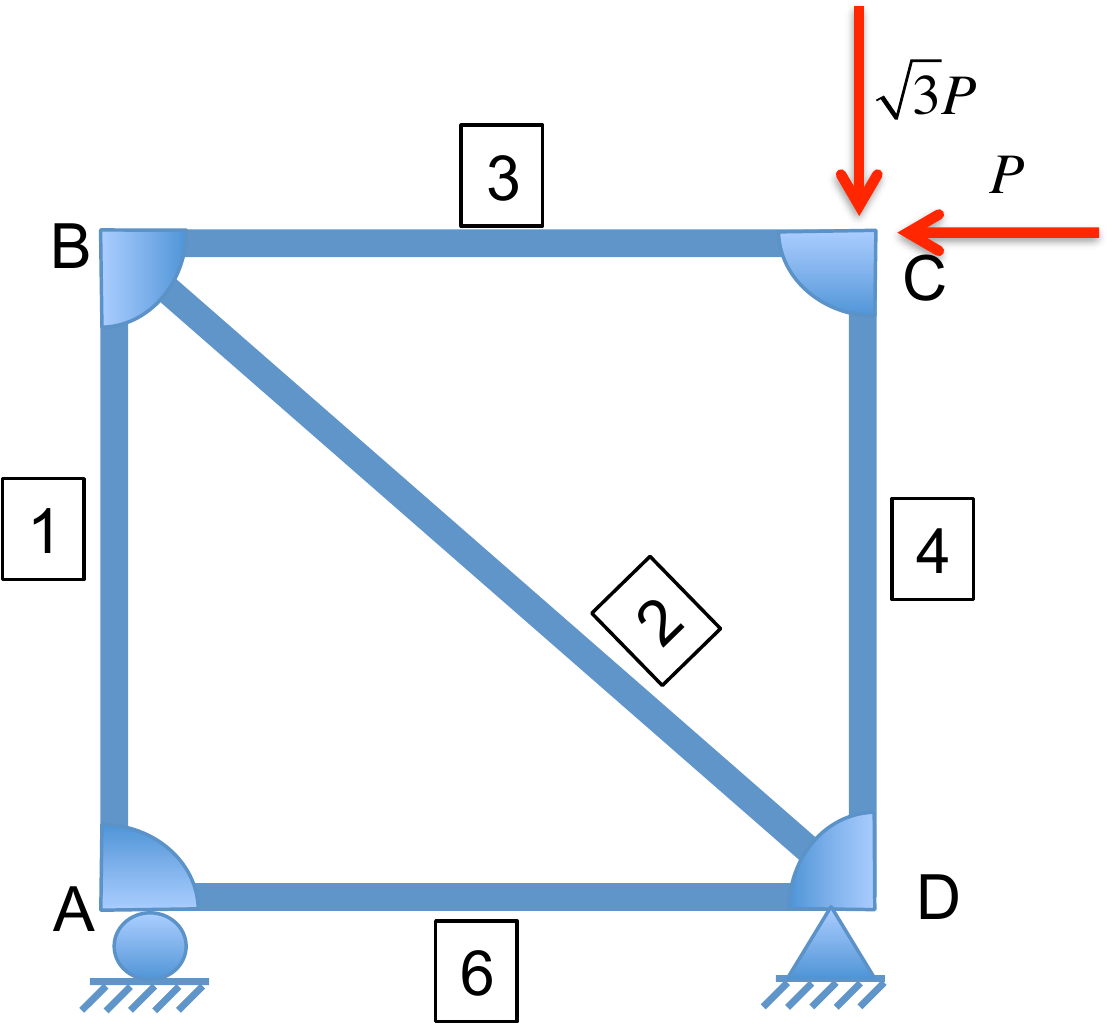}
\caption{Statically determinate}
\label{fig:subim341}
\end{subfigure}
\caption{Truss structures with and without member 5}
\label{fig:image4}
\end{figure}
\subsubsection{Two Dimensional random input}
In this section, the diagonal members 2 and 5 are assumed to have random cross-sectional areas, subject to the uniform distributions $A_2$ $\sim$ U(3, 9) $cm^2$ and $A_5$ $\sim$ U(3, 9) $cm^2$. The input parameters given in Table 1.
\begin{table}[h!]
  \centering
  \caption{Parameter values for the 2D truss problem}
  \label{tab:table1}
  \begin{tabular}{cccc}
    \toprule
    Members & Young's Modulus(E,GPa) & Area (A,$cm^2$) & Length(m)\\
    \midrule
    1 & $200$ & 6 & $2\sqrt{3}$\\
    2 & $200$  & U(3,9) & $4$\\
    3 & $200$  & 6 & $2$\\
    4 & $200$  & 6 & $2\sqrt{3}$\\
    5 & $200$  & U(3,9) & $4$\\
    6 & $200$  & 6 & $2$\\
    \bottomrule
  \end{tabular}
\end{table}

\begin{figure}[htbp]
\centering
\begin{subfigure}{0.4\textwidth}
\centering
\includegraphics[width=1.5\linewidth, height=7cm]{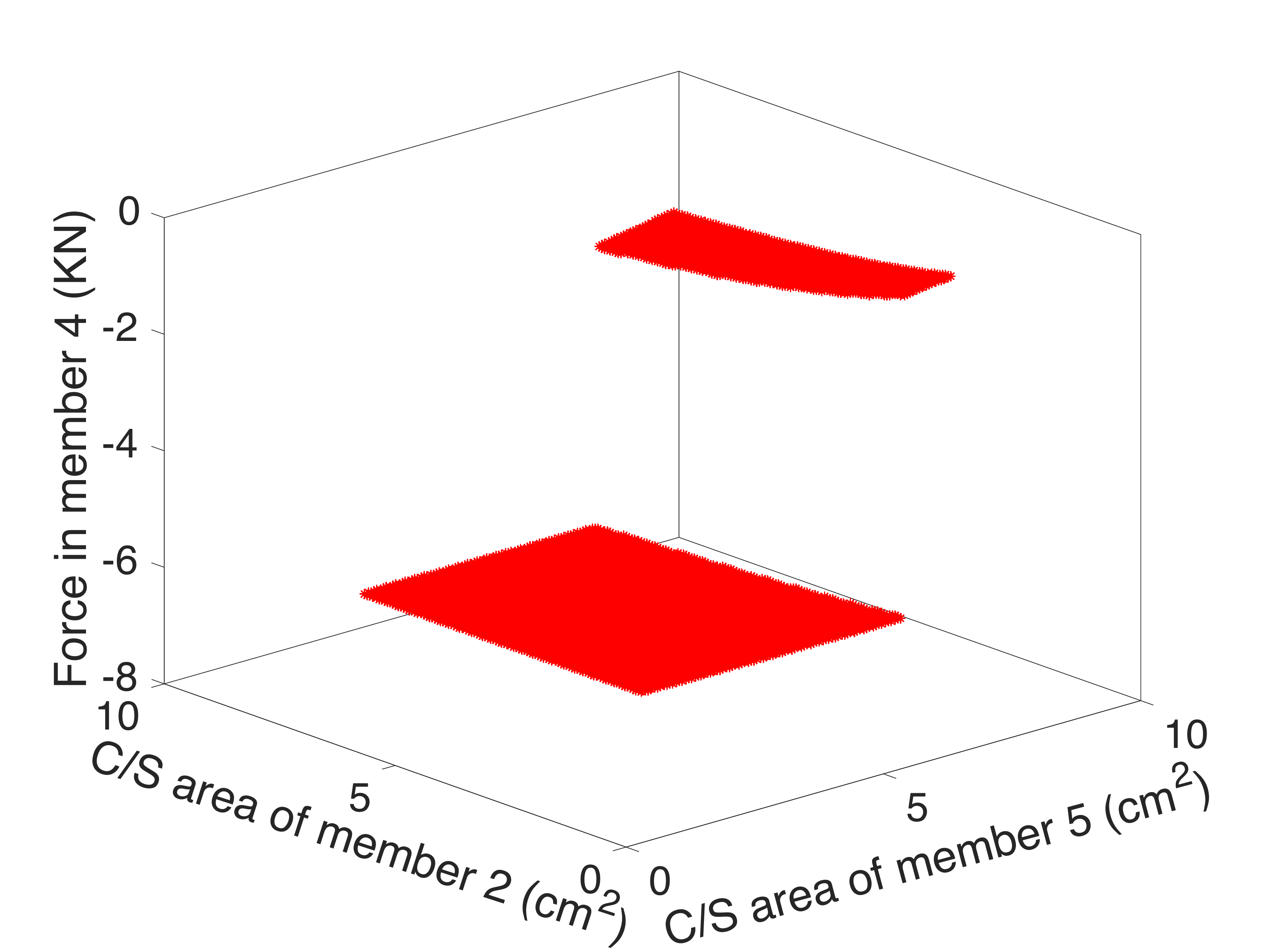} 
\caption{Exact function}
\label{fig:subim331}
\end{subfigure}
\hfill
\begin{subfigure}{0.4\textwidth}
\centering
\includegraphics[width=1.5\linewidth, height=7cm]{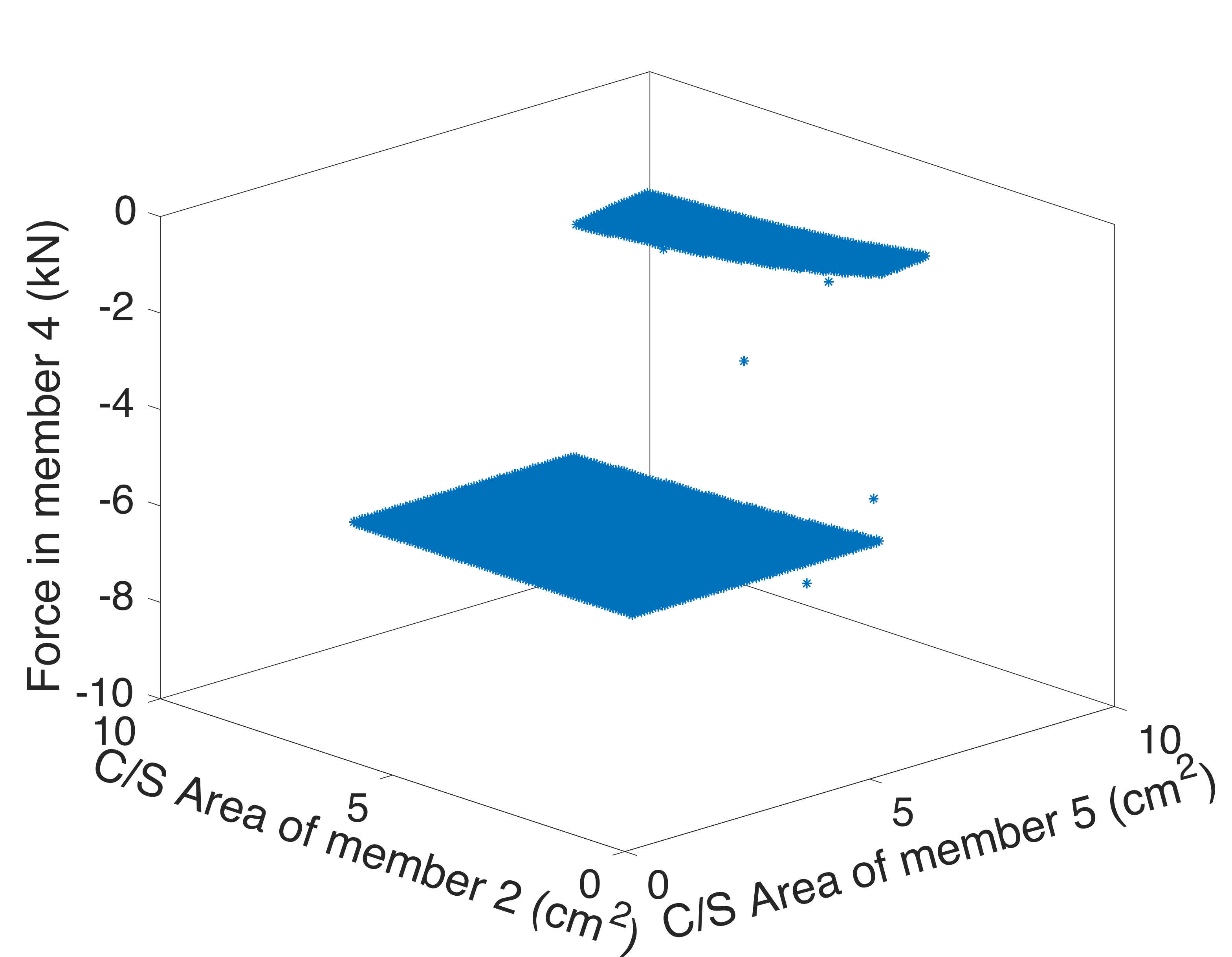}
\caption{Surrogate function}
\label{fig:subim341}
\end{subfigure}
\hfill
\begin{subfigure}{0.4\textwidth}
\centering
\includegraphics[width=1.5\linewidth, height=6cm]{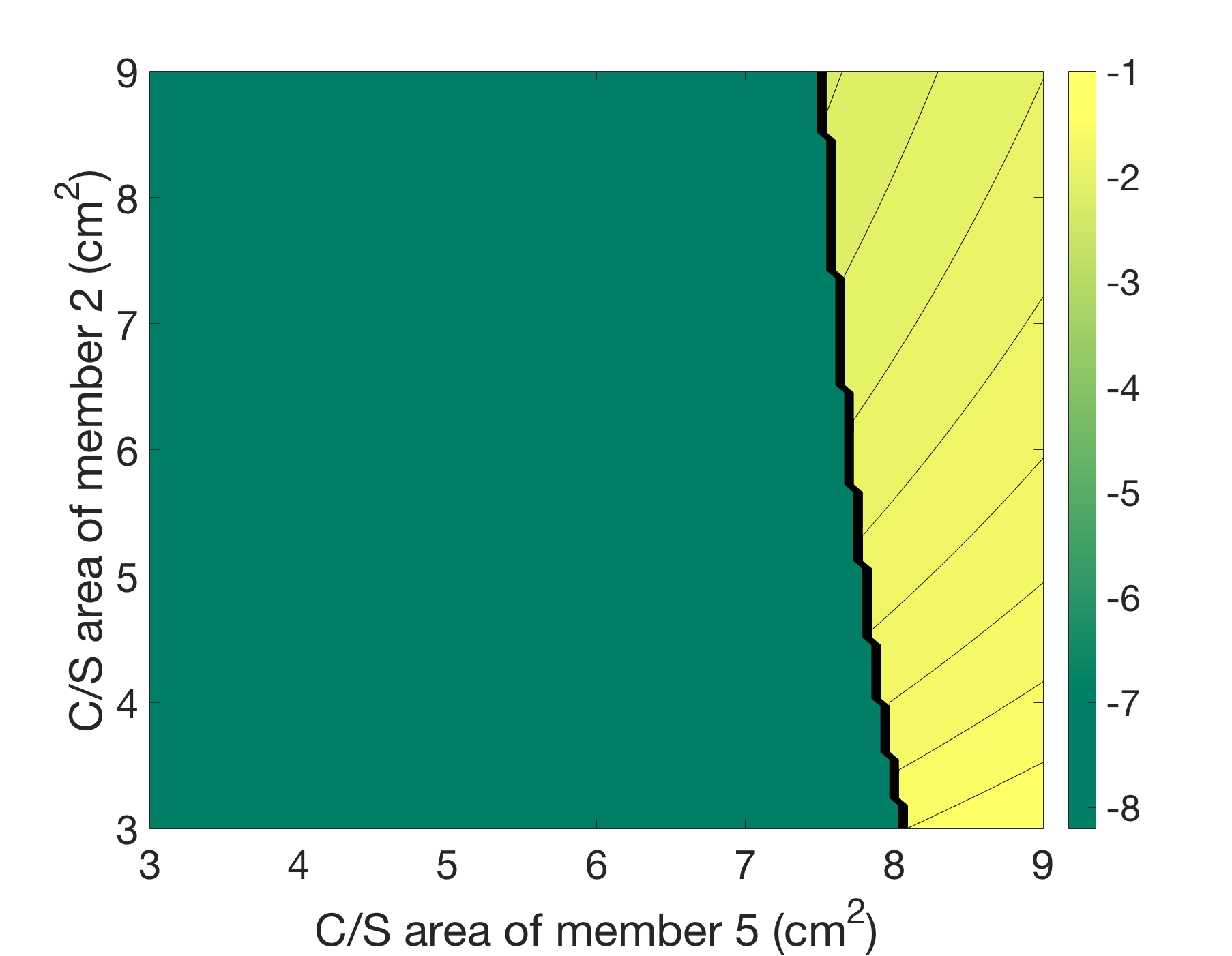}
\caption{Exact contour plot}
\label{fig:subim341}
\end{subfigure}
\hfill
\begin{subfigure}{0.4\textwidth}
\centering
\includegraphics[width=1.5\linewidth, height=6cm]{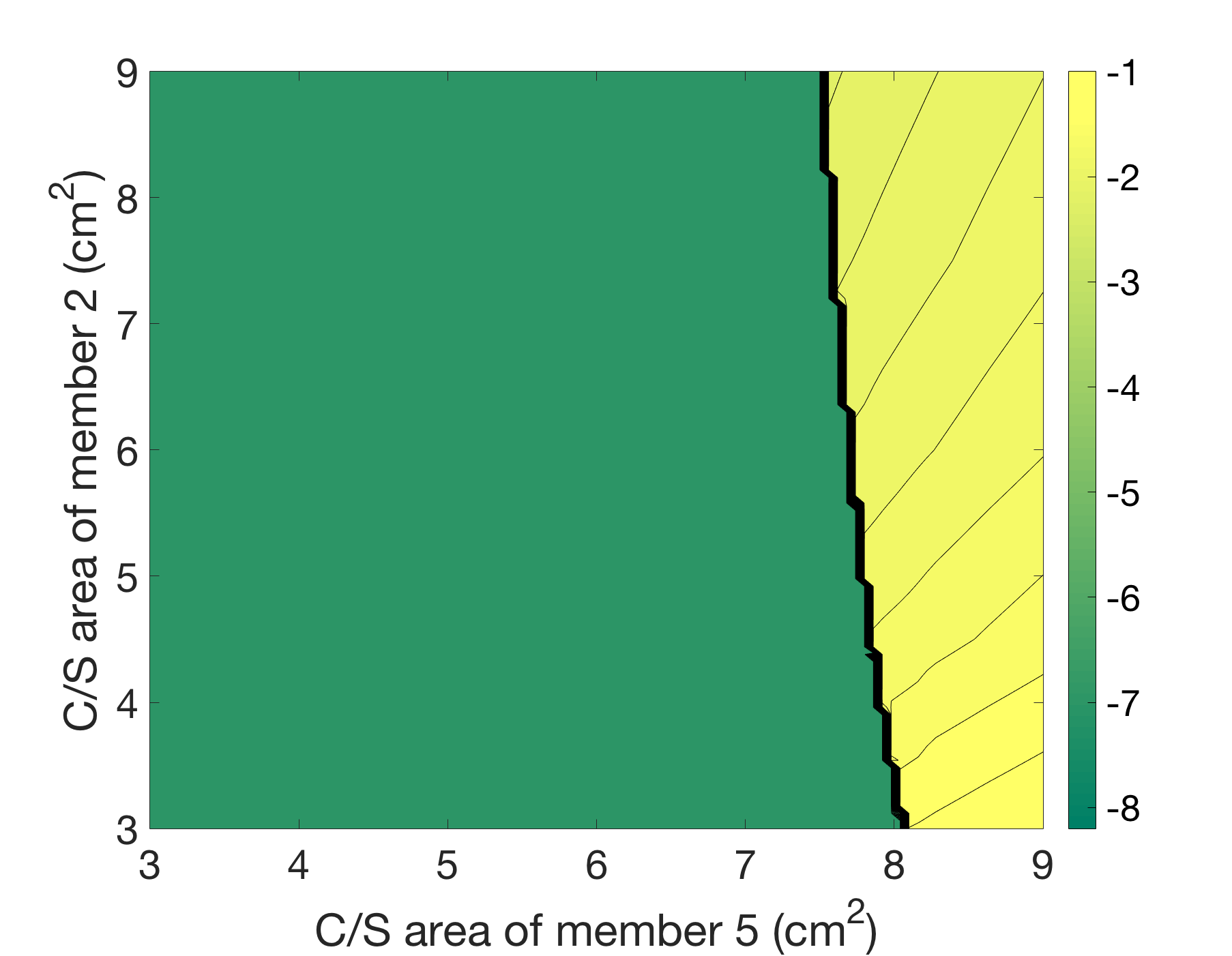}
\caption{Surrogate contour plot}
\label{fig:subim341}
\end{subfigure}
\caption{2 dimensional truss problem output}
\label{fig:image4}
\end{figure}
The true variation in the output of interest is shown in Figure 8(a). The surrogate response surface generated by the E-ASGC approach is shown in Figure 8(b). The corresponding contour plots are shown in Figure 8(c) and 8(d) respectively. Figure 9(a) shows how the input domain is sampled using the E-ASGC method. The E-ASGC method is compared with the ASGC method in the convergence plot of the root mean square error shown in Figure 9(b). The E-ASGC method clearly outperforms the ASGC method especially when the accuracy level is higher. For a root mean square error of 0.0053, the E-ASGC method requires 6,135 function evaluations while the ASGC method requires 12,008 function evaluations, thus reducing the sampling by almost a relative factor of 2.\\

\begin{figure}
\centering
\begin{subfigure}{0.4\textwidth}
\centering
\includegraphics[width=1.5\linewidth, height=6cm]{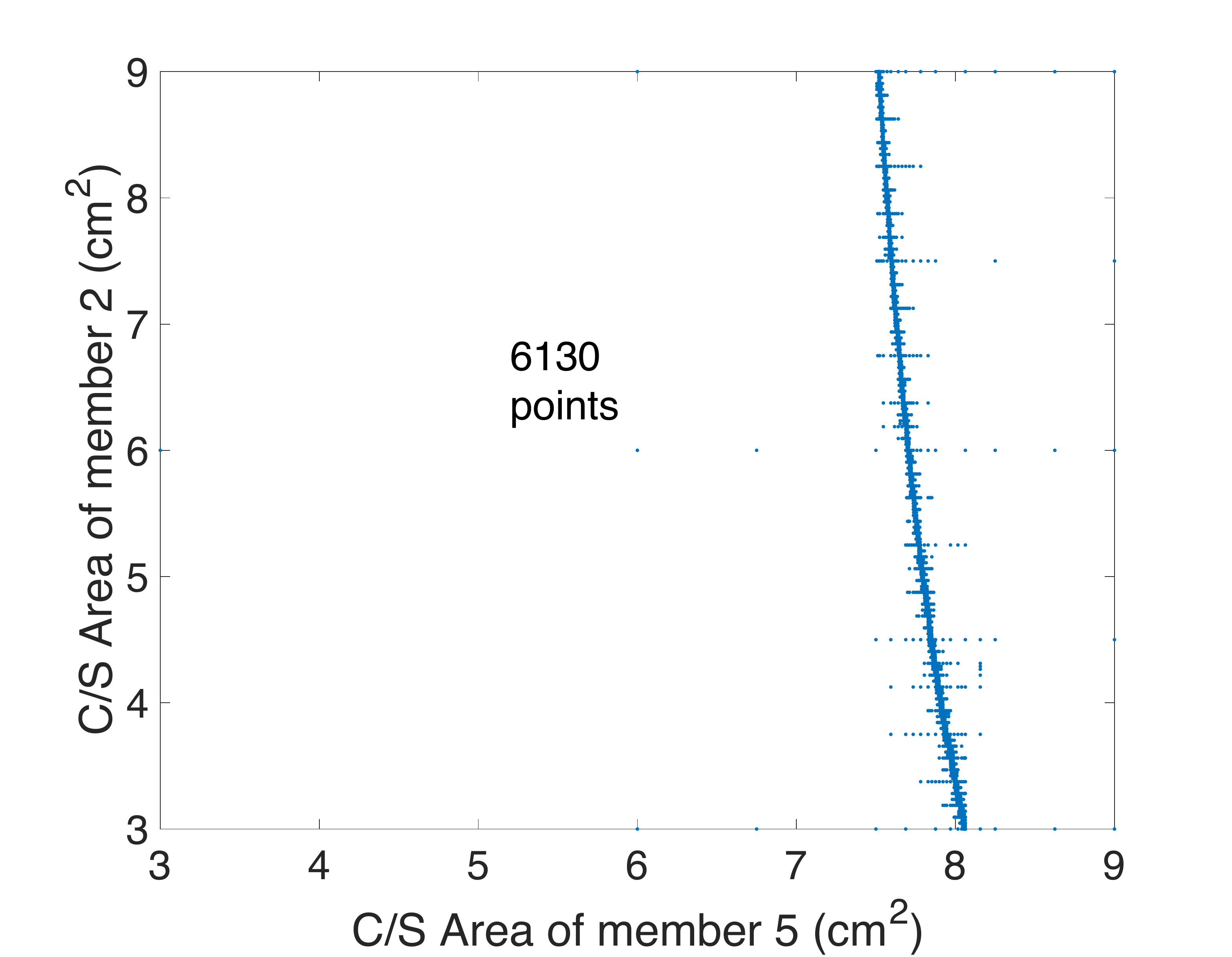}
\caption{E-ASGC input domain}
\label{fig:subim341}
\end{subfigure}
\hfill
\begin{subfigure}{0.4\textwidth}
\centering
\includegraphics[width=1.5\linewidth, height=6cm]{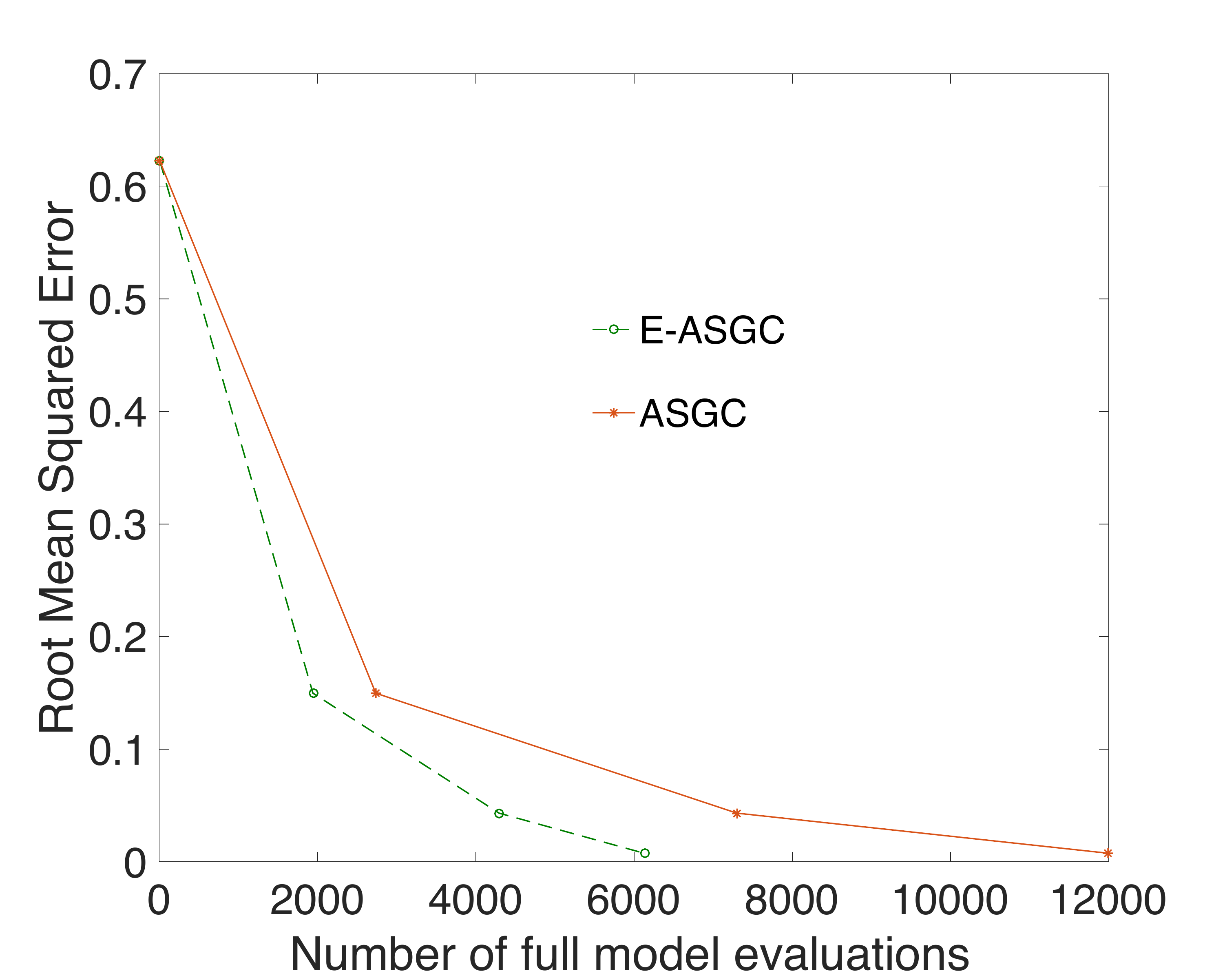}
\caption{Converge plot comparison}
\label{fig:subim341}
\end{subfigure}
\caption{Input Domain and Convergence plot for 2 dimensional truss problem}
\label{fig:image5}
\end{figure}
\subsubsection{3 Dimensional random input}
In this section, truss members 1, 3 and 5 are assumed to have random cross-sectional areas. The cross-section of the diagonal elements 5 is subject to the uniform distribution $A_5$ $\sim$ U(3, 9) $cm^2$ while the horizontal and vertical members 1 and 3 have cross-sectional areas all subject to uniform distributions $A_1$ $\sim$ U(5.5,6.5) $cm^2$ and $A_3$ $\sim$ U(5.5,6.5) $cm^2$ respectively. The input parameters are given in Table 2.
\begin{table}[h!]
  \centering
  \caption{Parameter values for the 3D truss problem}
  \label{tab:table1}
  \begin{tabular}{cccc}
    \toprule
    Members & Young's Modulus(E,GPa) & Area (A,$cm^2$) & Length(m)\\
    \midrule
    1 & $200$  & U(5.5,6.5) & $2\sqrt{3}$\\
    2 & $200$  & 6 & $4$\\
    3 & $200$  & U(5.5,6.5) & $2$\\
    4 & $200$  & 6 & $2\sqrt{3}$\\
    5 & $200$  & U(3,9)& $4$\\
    6 & $200$  & 6 & $2$\\
    \bottomrule
  \end{tabular}
\end{table}

Here, the E-ASGC method is again compared with the ASGC method in the convergence plot of the root mean square error shown in Figure 10. It is seen that with increase in accuracy, the efficiency of the E-ASGC method relative to the ASGC method also increases. For a root mean square error of 0.0580, the E-ASGC method requires 47,364 function evaluations while the ASGC method requires 154,677 function evaluations, thus reducing the sampling by a relative factor of around 3.\\


\begin{figure}[htbp]
\centering
\includegraphics[scale=0.3]{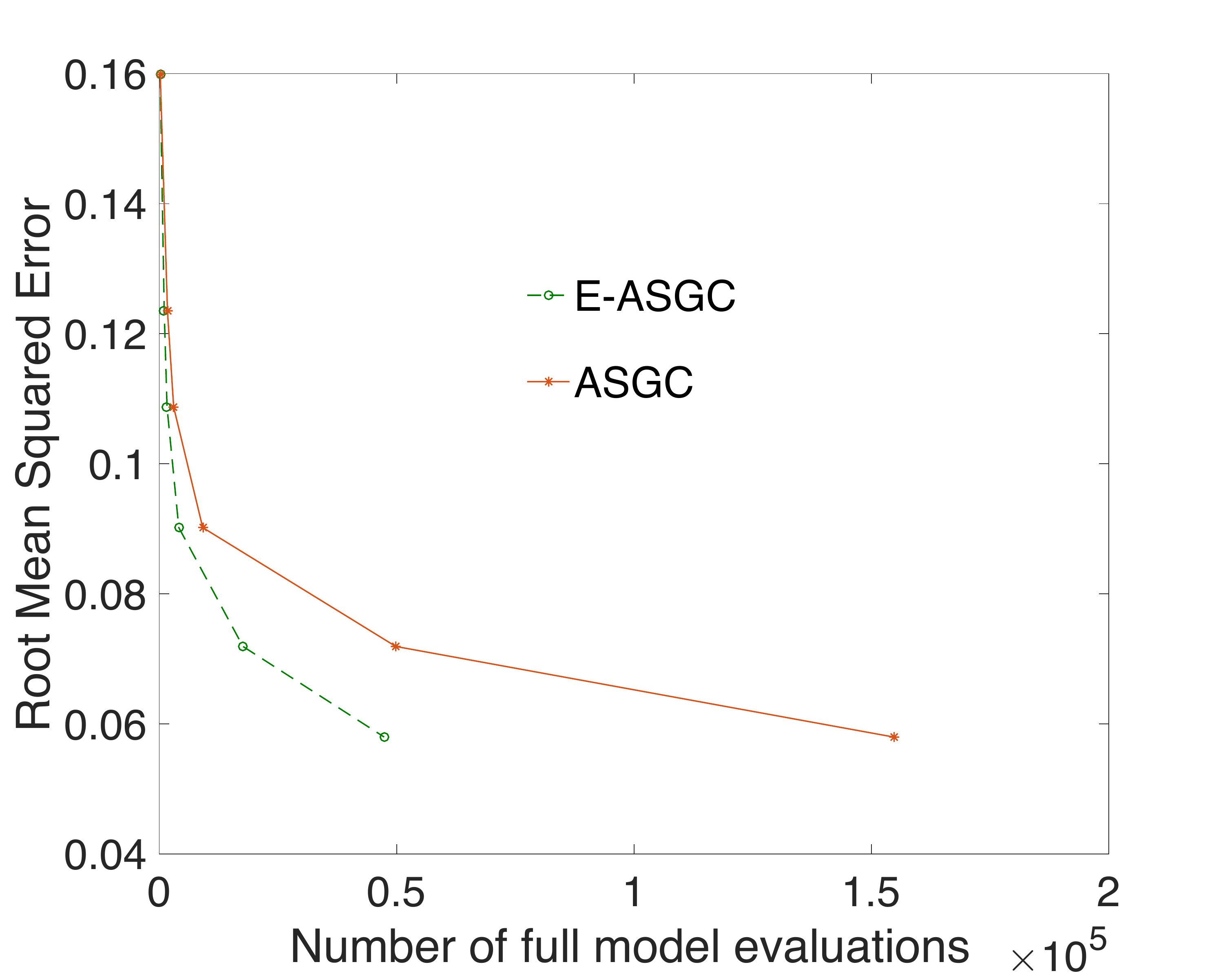}
\caption{Root Mean Squared Error convergence plot for 3 dimensional truss problem}
\end{figure}

\section{Conclusions} An efficient adaptive sparse grid approach through derivative estimation is developed which is based on the adaptive sparse grid subset collocation method (ASGC) \cite{ma2009adaptive}, which achieves faster convergence in the case of response functions that exhibit highly localized variations (such as discontinuities) in some regions and gradual variations in other regions of the stochastic input domain. The approach is significantly more efficient than the conventional sparse grid approaches. It is at least as efficient as the adaptive sparse grid subset collocation approach (ASGC), with up to two-fold increase in efficiency, depending on the nature of the response surface. If the response surface has many sharp variations, then the 1-D cubic spline interpolation cannot be suitably implemented and the efficiency reduces to that of the adaptive sparse grid subset approach. Also with increase in dimensions, the effectiveness relative to the ASGC method may decrease as the polynomial interpolations cover less space in the high-dimensional domain, given the control parameters remain the same as in the lower dimensional case. However, it is worth mentioning that given a complex deterministic model, any reduction in full model evaluations can be a significant contribution towards reducing computational cost.

\indent Future efforts will try to identify more efficient algorithms to segregate the smoother regions of the stochastic space from steep and discontinuous zones. Also, different interpolation techniques will be assessed that may achieve more accuracy in the smoother regions.

\section*{Acknowledgements} Research was sponsored by the Army Research Laboratory and was accomplished under
Cooperative Agreement Number W911NF-12-2-0022. The views and conclusions contained
in this document are those of the authors and should not be interpreted as representing the official
policies, either expressed or implied, of the Army Research Laboratory or the U.S. Government. The
U.S. Government is authorized to reproduce and distribute reprints for Government purposes
notwithstanding any copyright notation herein.\\
\indent The authors would like to thank Michael D. Shields, Robert M. Kirby and Yanyan He for their valuable suggestions.

\bibliographystyle{ieeetr}
\bibliography{reference}

\begin{thebibliography}{10}
\expandafter\ifx\csname url\endcsname\relax
  \def\url#1{\texttt{#1}}\fi
\expandafter\ifx\csname urlprefix\endcsname\relax\def\urlprefix{URL }\fi
\expandafter\ifx\csname href\endcsname\relax
  \def\href#1#2{#2} \def\path#1{#1}\fi

\bibitem{le2010spectral}
O.~Le~Ma{\^\i}tre, O.~M. Knio, Spectral methods for uncertainty quantification:
  with applications to computational fluid dynamics, Springer Science \&
  Business Media, 2010.

\bibitem{ma2009adaptive}
X.~Ma, N.~Zabaras, An adaptive hierarchical sparse grid collocation algorithm
  for the solution of stochastic differential equations, Journal of
  Computational Physics 228~(8) (2009) 3084--3113.

\bibitem{helton2003latin}
J.~C. Helton, F.~J. Davis, Latin hypercube sampling and the propagation of
  uncertainty in analyses of complex systems, Reliability Engineering \& System
  Safety 81~(1) (2003) 23--69.

\bibitem{rubinstein2011simulation}
R.~Y. Rubinstein, D.~P. Kroese, Simulation and the Monte Carlo method, Vol.
  707, John Wiley \& Sons, 2011.

\bibitem{au2003important}
S.~Au, J.~Beck, Important sampling in high dimensions, Structural safety 25~(2)
  (2003) 139--163.

\bibitem{niederreiter1992random}
H.~Niederreiter, Random number generation and quasi-monte carlo methods,
  cbms-nsf reg, in: Conf. Series Appl. Math, Vol.~63, 1992.

\bibitem{tu1996advantages}
J.~V. Tu, Advantages and disadvantages of using artificial neural networks
  versus logistic regression for predicting medical outcomes, Journal of
  clinical epidemiology 49~(11) (1996) 1225--1231.

\bibitem{pope1997computationally}
S.~B. Pope, Computationally efficient implementation of combustion chemistry
  using in situ adaptive tabulation.

\bibitem{shepard1968two}
D.~Shepard, A two-dimensional interpolation function for irregularly-spaced
  data, in: Proceedings of the 1968 23rd ACM national conference, ACM, 1968,
  pp. 517--524.

\bibitem{wiener1938homogeneous}
N.~Wiener, The homogeneous chaos, American Journal of Mathematics 60~(4) (1938)
  897--936.

\bibitem{ghanem2003stochastic}
R.~G. Ghanem, P.~D. Spanos, Stochastic finite elements: a spectral approach,
  Courier Corporation, 2003.

\bibitem{ghanem1998scales}
R.~Ghanem, Scales of fluctuation and the propagation of uncertainty in random
  porous media, Water Resources Research 34~(9) (1998) 2123--2136.

\bibitem{ghanem1999stochastic}
R.~Ghanem, Stochastic finite elements with multiple random non-gaussian
  properties, Journal of Engineering Mechanics 125~(1) (1999) 26--40.

\bibitem{ghosh2008strain}
D.~Ghosh, C.~Farhat, Strain and stress computations in stochastic finite
  element methods, International Journal for Numerical Methods in Engineering
  74~(8) (2008) 1219--1239.

\bibitem{doostan2009least}
A.~Doostan, G.~Iaccarino, A least-squares approximation of partial differential
  equations with high-dimensional random inputs, Journal of Computational
  Physics 228~(12) (2009) 4332--4345.

\bibitem{ghanem1999propagation}
R.~Ghanem, J.~Red-Horse, Propagation of probabilistic uncertainty in complex
  physical systems using a stochastic finite element approach, Physica D:
  Nonlinear Phenomena 133~(1) (1999) 137--144.

\bibitem{xiu2002wiener}
D.~Xiu, G.~E. Karniadakis, The wiener--askey polynomial chaos for stochastic
  differential equations, SIAM journal on scientific computing 24~(2) (2002)
  619--644.

\bibitem{xiu2002modeling}
D.~Xiu, G.~E. Karniadakis, Modeling uncertainty in steady state diffusion
  problems via generalized polynomial chaos, Computer methods in applied
  mechanics and engineering 191~(43) (2002) 4927--4948.

\bibitem{xiu2003modeling}
D.~Xiu, G.~E. Karniadakis, Modeling uncertainty in flow simulations via
  generalized polynomial chaos, Journal of computational physics 187~(1) (2003)
  137--167.

\bibitem{xiu2003new}
D.~Xiu, G.~E. Karniadakis, A new stochastic approach to transient heat
  conduction modeling with uncertainty, International Journal of Heat and Mass
  Transfer 46~(24) (2003) 4681--4693.

\bibitem{kewlani2009multi}
G.~Kewlani, K.~Iagnemma, A multi-element generalized polynomial chaos approach
  to analysis of mobile robot dynamics under uncertainty, in: 2009 IEEE/RSJ
  International Conference on Intelligent Robots and Systems, IEEE, 2009, pp.
  1177--1182.

\bibitem{wan2005adaptive}
X.~Wan, G.~E. Karniadakis, An adaptive multi-element generalized polynomial
  chaos method for stochastic differential equations, Journal of Computational
  Physics 209~(2) (2005) 617--642.

\bibitem{wan2006multi}
X.~Wan, G.~E. Karniadakis, Multi-element generalized polynomial chaos for
  arbitrary probability measures, SIAM Journal on Scientific Computing 28~(3)
  (2006) 901--928.

\bibitem{babuska2004galerkin}
I.~Babuska, R.~Tempone, G.~E. Zouraris, Galerkin finite element approximations
  of stochastic elliptic partial differential equations, SIAM Journal on
  Numerical Analysis 42~(2) (2004) 800--825.

\bibitem{le2004uncertainty}
O.~Le~Ma{\i}tre, O.~Knio, H.~Najm, R.~Ghanem, Uncertainty propagation using
  wiener--haar expansions, Journal of Computational Physics 197~(1) (2004)
  28--57.

\bibitem{ghosh2008stochastic}
D.~Ghosh, R.~Ghanem, Stochastic convergence acceleration through basis
  enrichment of polynomial chaos expansions, International journal for
  numerical methods in engineering 73~(2) (2008) 162--184.

\bibitem{xiu2010numerical}
D.~Xiu, Numerical methods for stochastic computations: a spectral method
  approach, Princeton University Press, 2010.

\bibitem{babuvska2007stochastic}
I.~Babu{\v{s}}ka, F.~Nobile, R.~Tempone, A stochastic collocation method for
  elliptic partial differential equations with random input data, SIAM Journal
  on Numerical Analysis 45~(3) (2007) 1005--1034.

\bibitem{mathelin2003stochastic}
L.~Mathelin, M.~Y. Hussaini, T.~A. Zang, A stochastic collocation algorithm for
  uncertainty analysis.

\bibitem{cools2002advances}
R.~Cools, Advances in multidimensional integration, Journal of Computational
  and Applied Mathematics 149~(1) (2002) 1--12.

\bibitem{bungartz2004sparse}
H.-J. Bungartz, M.~Griebel, Sparse grids, Acta numerica 13 (2004) 147--269.

\bibitem{smolyak1963quadrature}
S.~A. Smolyak, Quadrature and interpolation formulas for tensor products of
  certain classes of functions, in: Dokl. Akad. Nauk SSSR, Vol.~4, 1963, p.
  123.

\bibitem{klimke2005algorithm}
A.~Klimke, B.~Wohlmuth, Algorithm 847: spinterp: Piecewise multilinear
  hierarchical sparse grid interpolation in matlab, ACM Transactions on
  Mathematical Software (TOMS) 31~(4) (2005) 561--579.

\bibitem{klimke2006uncertainty}
W.~A. Klimke, Uncertainty modeling using fuzzy arithmetic and sparse grids,
  Citeseer, 2006.

\bibitem{klimke2007sparse}
A.~Klimke, Sparse grid interpolation toolbox--user’s guide, IANS report 17.

\bibitem{ganapathysubramanian2007sparse}
B.~Ganapathysubramanian, N.~Zabaras, Sparse grid collocation schemes for
  stochastic natural convection problems, Journal of Computational Physics
  225~(1) (2007) 652--685.

\bibitem{nobile2008anisotropic}
F.~Nobile, R.~Tempone, C.~G. Webster, An anisotropic sparse grid stochastic
  collocation method for partial differential equations with random input data,
  SIAM Journal on Numerical Analysis 46~(5) (2008) 2411--2442.

\bibitem{foo2008multi}
J.~Foo, X.~Wan, G.~E. Karniadakis, The multi-element probabilistic collocation
  method (me-pcm): Error analysis and applications, Journal of Computational
  Physics 227~(22) (2008) 9572--9595.

\bibitem{agarwal2009domain}
N.~Agarwal, N.~R. Aluru, A domain adaptive stochastic collocation approach for
  analysis of mems under uncertainties, Journal of Computational Physics
  228~(20) (2009) 7662--7688.

\bibitem{chantrasmi2009pade}
T.~Chantrasmi, A.~Doostan, G.~Iaccarino, Pad{\'e}--legendre approximants for
  uncertainty analysis with discontinuous response surfaces, Journal of
  Computational Physics 228~(19) (2009) 7159--7180.

\bibitem{elman2012stochastic}
H.~C. Elman, C.~W. Miller, Stochastic collocation with kernel density
  estimation, Computer Methods in Applied Mechanics and Engineering 245 (2012)
  36--46.

\bibitem{grigoriu2014response}
M.~Grigoriu, Response statistics for random heterogeneous microstructures,
  SIAM/ASA Journal on Uncertainty Quantification 2~(1) (2014) 252--275.

\bibitem{kincaid2002numerical}
D.~R. Kincaid, E.~W. Cheney, Numerical analysis: mathematics of scientific
  computing, Vol.~2, American Mathematical Soc., 2002.

\bibitem{hall1976optimal}
C.~A. Hall, W.~W. Meyer, Optimal error bounds for cubic spline interpolation,
  Journal of Approximation Theory 16~(2) (1976) 105--122.

\bibitem{genz1987package}
A.~Genz, A package for testing multiple integration subroutines, in: Numerical
  Integration, Springer, 1987, pp. 337--340.

\end{thebibliography}

\end{document}